\documentclass[a4paper,11pt, reqno]{amsart}

\usepackage[ansinew]{inputenc}
\usepackage{amsmath,amsthm,amsfonts,amssymb,stackengine}
\usepackage{mathrsfs}
\usepackage[left=2.50cm,right=2.25cm,top=2.50cm,bottom=2.50cm]{geometry}
\usepackage{graphicx}
\usepackage[english]{babel}
\usepackage{enumitem,fancyhdr}
\usepackage{etoolbox,lineno}

\usepackage{fouridx}
\usepackage{bigints}
\usepackage{soul}
\usepackage{bbm}

\usepackage[dvipsnames]{xcolor}

\newcommand\red[1]{{\color{red}#1}}

\numberwithin{equation}{section}
\theoremstyle{plain}
\newtheorem{Theorem}{Theorem}[section]

\newtheorem{lemma}[Theorem]{Lemma}

\theoremstyle{definition}
\newtheorem{definition}[Theorem]{Definition}
\newtheorem{remark}[Theorem]{Remark}

\newcommand\del[1]{}

\del{

	\newcommand\del[1]{}
	}

\DeclareOldFontCommand{\rm}{\normalfont\rmfamily}{\mathrm}

\newcommand{\CP}{{{ \mathcal P }}}

\newcommand{\CF}{{{ \mathcal F }}}

\newcommand{\lk}{\left}
\newcommand{\lqq}{\lefteqn}
\newcommand{\rk}{\right}

\newcommand{\DEQS}{\begin{eqnarray*}}
	\newcommand{\EEQS}{\end{eqnarray*}}
\newcommand{\DEQSZ}{\begin{eqnarray}}
	\newcommand{\EEQSZ}{\end{eqnarray}}
\newcommand{\DEQ}{\begin{eqnarray}}
	\newcommand{\EEQ}{\end{eqnarray}}

\allowdisplaybreaks
\allowdisplaybreaks
\def\XXint#1#2#3{{\setbox0=\hbox{$#1{#2#3}{\int}$ }
		\vcenter{\hbox{$#2#3$ }}\kern-.6\wd0}}

\usepackage{parskip}

\usepackage{scalerel,stackengine}
\stackMath
\newcommand\reallywidehat[1]{%
	\savestack{\tmpbox}{\stretchto{%
			\scaleto{%
				\scalerel*[\widthof{\ensuremath{#1}}]{\kern-.6pt\bigwedge\kern-.6pt}%
				{\rule[-\textheight/2]{1ex}{\textheight}}
			}{\textheight}%
		}{0.5ex}}%
	\stackon[1pt]{#1}{\tmpbox}%
}

\begin{document}
	 \title{ {On the inertial range bounds of K-41 like Magnetohydrodynamics turbulence.}}
	 \author{T.A. Tegegn}
	 \address{T.A. Tegegn, Department of Mathematics and Applied Mathematics,
	 	Sefako Makgatho Health Sciences University, South Africa.}
	 \email{tesfalem.tegegn@smu.ac.za}
	 \maketitle
	 
	 \begin{abstract}
	 	 {The spectral slope of Magnetohydrodynamic (MHD) turbulence varies depending on the spectral theory considered; $ -3/2 $ is the spectral slope in Kraichnan-Iroshnikov-Dobrowolny (KID)   theory, $ -5/3 $ in Marsch-Matthaeus-Zhou's  and Goldreich-Sridhar theories also called Kolmogorov-like  (K-41 like) MHD theory, combination of the $-5/3$ and $-3/2$ scales in Biskamp and so on. A rigorous mathematical proof to any of these spectral theories is of great scientific interest. Motivated by the 2012 work of A. Biryuk and W. Craig [Physica D 241(2012) 426-438], we  establish inertial range bounds for K-41 like   phenomenon in MHD turbulent flow through a mathematical rigour;  a range of wave numbers in which the spectral slope of MHD turbulence is proportional to $-5/3$ is established and the upper and lower bounds of this range are explicitly formulated. We also have shown that the Leray weak solution of the standard MHD model is bonded in the Fourier space,  the spectral energy of the system is bounded and its average over time decreases in time.  }
	 \end{abstract}
 
 \textbf{Keywords:} {Magnetohydrodynamics turbulence; Harmonic  Analysis; Kolmogorov theory; Inertial range bound
 
 
 \vspace*{1cm}
 
 \hrule
 
 \vspace*{1cm}
 
 
 \pagestyle{fancy}
 \fancyhead{}\fancyfoot{}

 \lhead{T.A. Tegegn}
 \rhead{Inertial Bounds in  MHD Turbulence}
 \cfoot{\arabic{page}}

 \tableofcontents
 \section{Introduction}
 \label{introduction}
 
 At high Reynolds number fluid  and plasma flows exhibit a complex random behavior  called turbulence. 
 Turbulence is observed in a great majority of 
 fluids both in nature such as the atmosphere, river currents, oceans, solar wind and interstitial bodies and in technical devices, such as laboratory installations, nuclear power plants, etc. 
 Its importance in industry and physical sciences, such as making predictions about heat transfer in nuclear power plants, drag in oil pipelines and the weather is tremendous. 
 Besides these real life relevant issues the study of turbulence can assist mathematical researchers in understanding some aspects, such as regularity of Euler's equation, 
 Navier-Stokes equation,  Magnetohydrodynamics equations and so on, see for instance \cite{chen2012kolmogorov}.

 \del{\red{According to literature, many generation of scientists have been studying turbulence to unlock its mysteries   ever since the very systematic observation by Leonardo da Vinci at the beginning of 16th century. }}
  {Literature shows that the phenomenon of turbulence has captured the attention of humankind for centuries, see for instance \cite{davidson2011voyage}. 
 	The discovery of Euler equations in the mid of the 18th century   
 	and Navier-Stokes equations in the first half of the 19th century  
 	are the major scientific and mathematical breakthrough developments in terms of having  governing rule for   fluid flows.
 	Towards the end of 19th century Osborne Reynolds  
 	laid a foundation for the theory of turbulence, 
 	see  \cite{jackson2007osborne}\cite{reynolds1883experimental,reynolds1894dynamical}\cite[p.~488]{tikhomirov1991selected}.} 
 Reynold's number,  a widely used criteria  to classify whether a given flow is turbulent or not, and 
 Reynolds averaged Navier-Stokes  equations (RANS) are due to O. Reynolds.
 \del{see \cite{reynolds1883experimental,reynolds1894dynamical}\cite[p.~488]{tikhomirov1991selected}\cite{jackson2007osborne}.
 	RANS is derived}
 RANS is formulated by decomposing the  velocity field $u(x,t)$ in to average velocity $\bar{u}(x,t)$ over a time interval and    fluctuation velocity $u'(x,t) = u(x,t)-\bar{u}(x,t)$, 
 and finally rewriting Navier-Stokes equations in terms of the average velocity $\bar{u}$. In fact, RANS is still one of the most widely used models to study turbulence in fluids, see \cite{alfonsi2009reynolds,argyropoulos2015recent} and the references there.

 \del{\red{  {MOSTLY TO BE REMOVED} In fact, several great contribution toward the development of turbulence theories   have been made by a series of scientists, such as Ludwig Prandtl \cite{Prandtl1904}, Theodore von K\'{a}rm\'{a}n \cite{von1948progress,von1949concept}, Geoffrey Ingram Taylor \cite{taylor1935statistical,taylor1921experiments}, Millionshchikov \cite{millionshchikov1939decay}, Andrej Nikolaevich Kolmogorov \cite{kolmogorov1941degeneration,kolmogorov1941dissipation,kolmogorov1941equations,kolmogorov1941local} and many more. 
 		The book ``A Voyage through turbulence'', \cite{davidson2011voyage} would be a nice reference in this regard.} }

 In 1941, in a series of works \cite{kolmogorov1941degeneration,kolmogorov1941dissipation,kolmogorov1941equations,kolmogorov1941local}, A.N. Kolmogorov set a phenomenological theory for hydrodynamic turbulence. His theory   postulated that the spectral energy of a fully developed turbulence decays according to the rule
 
 %
 %
 %
 \DEQSZ
 C_0\epsilon^{2/3}k^{-5/3} \label{mainspectra},
 \EEQSZ
 
 \noindent over a range of wave numbers, $k\in[k_1,k_2]$, also called the inertial range; where $\epsilon$ is the energy dissipation rate and $C_0$ is a universal constant called Kolmogorov constant. 
 The exponents in \eqref{mainspectra} are determined by dimensional analysis. 
 The theory is often referred to as K-41 theory or Kolmogorov's $-5/3$ law.
 %
 %
 The state of the art exposition of Komogorov's school of turbulence can be found in the seminal monographs of Monin and Yaglom \cite{moninyaglomI,monin2013statisticalII}.
 
 Works particularly focus on the spectral energy of MHD flows started to emerge by the mid of the 20th century.
 From the earliest of such works 
 the works of 
 Kraichnan \cite{kraichnan1965inertial,kraichnan1971inertial} and Iroshnikov \cite{iroshnikov1964turbulence} can be mentioned. 
 Unlike Kolmogorov where the spectral energy decays proportional to $k^{-5/3}$,  Kraichnan and Iroshnikov concluded that  
 the spectral energy  
 of a fully developed MHD turbulent flow decays proportional to $k^{-3/2}$, which  
 later  on was supported by M. Dobrowolny, A. Mangeney and P. Veltri in \cite{dobrowolny1980fully}. 
 Mahendra K. Verma  in his review  \cite{verma2004statistical} said that these works 
 are the first to establish  phenomenological theory on MHD turbulence,
 which he 
  {called  
 	Kraichnan-Iroshnikov-Dobrowolny} (KID) phenomenon.
 %
 %
 %
 %


  {
 	MHD turbulence unlike hydrodynamic turbulence is controlled by a combined effect of the magnetic field and the fluid velocity, see for instance \cite{Chandrasekhar1955}.
 	Despite the difference in the  formation of   hydrodynamic and MHD turbulence, several authors have argued that  
 	under certain conditions the spectral energy of MHD turbulence also decays proportional to $k^{-5/3},$ which is widely accepted as a spectral slope for hydrodynamic turbulence. For instance,   Marsch and Tu in \cite{marsch1990radial} and Marsch in \cite{marsch1991turbulence} suggested that the decay rate of an isotropic turbulence in solar wind is very likely to be $-5/3$ than $-3/2$.
 	 {Matthaeus and Zhou in  \cite{matthaeus1989extended}} proposed that the larger wave numbers (relative to the mean magnetic field) would follow the ${-3/2}$ law whereas the smaller wavenumbers  would follow the ${-5/3}$ law. Biskamp in \cite{biskamp1994cascade} proposed three different rates; ${-5/3}$ for the general MHD turbulence when   Alfv\'en effect is neglected,  $-5/4$ when   Alfv\'en effects are included and the mean magnetic field is constant   and $ {-3/2} $ when Alfv\'en effects are considered and the mean magnetic field is fluctuating. 
 	Boldyrev in \cite{boldyrev2005spectrum} also concluded that MHD turbulence is not completely described by either the $ -3/2 $ or $-5/3$ scales; the scales depend on the strength of the external magnetic field: $-3/2$ scale applies when the mean magnetic field is strong while $-5/3$ scale applies when the external magnetic field is weak.  
 	We refer to the review by Verma \cite{verma2004statistical} for the several phenomenological theories on MHD turbulence and 
 	the book by Davidson et. al. \cite{davidson2011voyage} 
 	for the biographies and works of some of the prominent contributors to the area. 
 	
 }
 

 \del{Despite the clear difference between  KID  and  K-41, 
 	some studies established that there are cases where MHD turbulence satisfies K-41. 
 	For instance, Verma in \cite{verma1999mean} has noted that his observational results and  calculations agreed with Kolmogorov's theory more than with the KID and Biskamp in \cite{biskamp1994cascade} 
 	concluded that in the general MHD case the behavior of  fully developed MHD turbulence  is close to K-41 rather than KID.}
 
 The purpose of  this paper is to  establish a spectral range for  K-41 like MHD phenomenon through mathematical rigour. 
 The work was motivated by the 2012 paper of Andei Biryuk and Walter Craig \cite{biryuk2012bounds} where they established
 an  estimate for the Leray  weak solution  of Navier-Stokes equations in the norm $ \|  \widehat{\partial_xu(\cdot,t)} \|_{L^\infty}$  which lead to prove  solution's ability to satisfy Kolmogorov's spectral law \eqref{mainspectra}. 
 Leray's  weak solution 
 was formulated in the first half of 1930's by J. Leray \cite{leray1931systeme,leray1933etude}.  
  {   
 	Interestingly, following Leray's work several authors treated   weak solutions for fluid dynamic models as  turbulent solution, see for instance   \cite{duchon2000inertial,eyink2001dissipation,filho2006weak,kato1984strongl,leray1931systeme,leray1933etude}. 
 	Therefore, it is not surprising to see the Leray weak solution of Navier-Stokes equations obeying K-41. In a similar passion, we consider the weak solution for a system of MHD equations   as a turbulent solution and attempt to show that it obeys the $-5/3$ spectral law over a range of wave numbers when certain condition are met. 
 }
 
 \del{Before we get in to the rigorous mathematical formulation, it is \red{worth noting} the main difference between hydrodynamic and MHD turbulence. It is obvious that in hydrodynamic turbulence  velocity or kinetic energy plays the central role but in MHD flows  turbulence is driven by  the combined effect of velocity and magnetic fields.
 	S. Chandrasekhar in \cite{Chandrasekhar1955} described the role played by the two fields as follows;
 	``$\ldots$the amplification of the magnetic field by the turbulent motions and the suppression of the motions by the magnetic field will balance each other and one may expect that an equipartition between the two forms of energy will result.''}
 
 The dynamics of MHD flows in general is described by a system of partial differential equations given by
 \DEQSZ  
 \lk\lbrace \begin{array}{l l}
 	\partial_t u +(u\cdot\nabla)u+\nabla\pi - (b\cdot \nabla)b-\nu\Delta u = f_1, & (0,\infty)\times D ,\\
 	\partial_t b +(u\cdot\nabla)b - (b\cdot \nabla)u-\eta\Delta b = f_2, & (0,\infty)\times D, \\
 	\text{div } u = {\rm div } b = 0 ,& D, \\
 	u|_{t=0} = u_0,  \quad b|_{t=0}=b_0,  & D,
 \end{array} \rk.   \label{mhdmain}
 \EEQSZ  
 
 \noindent where  $u=u(x,t)$ is the flow velocity, $b=b(x,t)$ is the magnetic field, 
  {$\pi=P+\frac12|b|^2$ is the total pressure on the system with $P$ representing the pressure function from the equation of motion}, $\nu>0$ is the kinetic viscosity of the fluid,  $\eta>0$ is the resistivity of the fluid and the spatial domain $D$ is  the Euclidean space $\mathbb{R}^3$. 
 The inhomogeneous external forces $f_1=f_1(x,t),\, f_2=f_2(x,t)$ are assumed to be divergence-free 
 and   satisfy $f_1,f_2\in L_{loc}^\infty([0,\infty);H^{-1}(D)\cap L^2(D))$,  {where $L^\infty_{loc}$ is the space of locally bounded functions, $H^{-1}$ and $L^{2}$ are the usual Sobolev and Lebesgue spaces, respectively}.
  {The derivation of equation \eqref{mhdmain} is done by combining the Navier-Stokes equations and the Maxwell equations in some way, see \cite{arsenio2014derivation,goedbloed1998derivation,peng2008rigorous}.
 }
 
 We now introduce 
 the spectral energy function, denoted by $E(k,t)$; the spectral energy of  the  MHD flow model \eqref{mhdmain} is 
 given by the surface integral
 
 \DEQSZ
 E(k,t) := \int_{|\xi|=k}(|\widehat{u}(\xi,t)|^2+|\widehat{b}(\xi,t)|^2) {\rm d} S(\xi), \quad k\in[0,\infty),\, \{|\xi| = k\} \subset D , \label{spectrafun}
 \EEQSZ
 where $\widehat{u}$ and $\widehat{b}$ represent the Fourier transforms of $u$ and $b$, respectively.

 Of great scientific interest is the question of rigorous  {mathematical} proof of the spectral theory, K-41 or otherwise, under physically admissible conditions. Therefore, our main goal will be  setting the conditions on the data and to show that  
 the spectral energy   \eqref{spectrafun} satisfies $-5/3$ law
  {when such    conditions are met.} 
 
 Before we give a formal definition to the weak solution of \eqref{mhdmain}, we  introduce some 
 function spaces and their notations as they appear in \cite{cannone2006cauchy}. 
 We denote by $C_{0,\sigma}^\infty$ the set of all divergence free smooth functions with compact support in $ D$.
 $L^p_\sigma$ is the closure of $ C_{0,\sigma}^\infty $ with respect to the $L^p$ norm in the usual sense. 
 For $1\leq p\leq\infty$ the space $L^p$ stands for the usual (vector-valued) Lebesgue space over $ \mathbb{R}^3.$
 For $s\in\mathbb{R}$, we denote by $H^s_\sigma$ the closure of $C_{0,\sigma}^\infty$ with respect to the $H^s$ norm. 
 \begin{definition}\label{mhdweaksol}
 	Let $(u_0,b_0)\in L_\sigma^2( D).$  A vector  $(u,b)$ is said to be a weak solution to \eqref{mhdmain} on $D\times[0,\infty) $ if it satisfies the following conditions:
 	\begin{enumerate}
 		\item   for any $T>0$ the vector function $(u,b)$ lies in the following function space,
 		
 		\DEQS
 		u,\, b\in L^\infty([0,T);L_\sigma^2(D))\cap  L^2([0,T);{H}_\sigma^1(D)).
 		\EEQS		
 		
 		\item  		the pair $(u,b)$ is a distributional solution of \eqref{mhdmain}; i.e., 
 		for every $(\Phi,\Psi)$ in  
 		
 		\DEQS
 		H^1((0,T); H^1_\sigma\cap L^2),
 		\EEQS
 		
 		with $\Phi(T)=\Psi(T)=0$,
 		
 		\DEQS
 		\int\limits_0^T&\{-(u,\partial_t\Phi) + \nu(\nabla u, \nabla\Phi) + (u\cdot\nabla u,\Phi)-(b\cdot \nabla b,\Phi) \}{\rm d} t\\& =-(u_0,\Phi(0)) + \int\limits_0^T (f_1,\Phi){\rm d} t ,
 		\EEQS
 		
 		and
 		
 		\DEQS
 		\int\limits_0^T&\{-(b,\partial_t\Psi)+\eta(\nabla b,\nabla\Psi) + (u\cdot\nabla b, \Psi)-(b\cdot\nabla u,\Psi) \}{\rm d} t \\
 		&= -(b_0,\Psi(0)) + \int\limits_0^T (f_2,\Psi){\rm d} t.
 		\EEQS
 		
 		Furthermore $\lim\limits_{t\rightarrow 0^+}u(\cdot,t) = u_0(\cdot)$ and $ \lim\limits_{t\rightarrow 0^+}b(\cdot,t) = b_0(\cdot)$ exist  in the strong $L^2 $ sense.
 		\item  
 		the following energy inequality is satisfied, 
 		
 		\DEQSZ
 		&&{\frac{1}{2}\int_D|u(x,t)|^2+|b(x,t)|^2{\rm d} x + \min({\nu,\eta})\int_0^t\int_D|\nabla u(x,s)|^2+|\nabla b(x,s)|^2{\rm d} x{\rm d} s} \nonumber\\
 		&&{-\int_0^t\int_Du(x,s)\cdot f_1(x,s)+b(x,s)\cdot f_2(x,s){\rm d} x{\rm d} s} \nonumber \\
 		&&\leq \frac{1}{2}\int_D|u_0(x)|^2+|b_0(x)|^2{\rm d} x \label{lerrayenergy}
 		\EEQSZ
 		
 	\end{enumerate}
 	for all $0<t<\infty$.
 \end{definition}

 The rest of the paper is organized as follows. In section \ref{sec: Fourier} we begin with a  brief  discussion on the Fourier transform and its properties   followed by  derivation of estimates for the weak solution of \eqref{mhdmain} in Fourier space. We present our main results together with their respective proofs in section \ref{sec: energyspectra}. We formulate a bound for the spectral energy function \eqref{spectrafun} and show that the time average of the bounds decreases in time. We also drive the spectral range bounds and formulate   necessary conditions on the data for the system \eqref{mhdmain} to exhibit K-41 like behaviour.

 \section{Estimates for the solution field $(u,b)$ in a Fourier Space}
 \label{sec: Fourier}

 \subsection{The Fourier transform}  
 
 The Fourier transform of an integrable function $u$, denoted    by $\widehat{u}$, 
 is defined   by
 
 \DEQS
 \widehat{u}(\xi) = \int_De^{-i\xi\cdot x}u(x){\rm d} x.
 \EEQS
 
 The Fourier transform has several interesting properties, among them the following three are of great importance to our work; 
 \DEQSZ
 \| u \|_{L^2(D)}^2 =  \| \widehat{u} \|_{L^2(D)}^2,\label{plancherel}
 \EEQSZ
 
  {
 	\DEQSZ
 	\widehat{  \partial_x^\alpha u(x) } = (i\xi)^\alpha \partial_\xi^\alpha \widehat{u}(\xi) \quad \text{and}\quad \widehat{x^\alpha u(x)} = (-i)^\alpha \partial_\xi^\alpha\widehat{u}(\xi),\label{Fourier-derivative}
 	\EEQSZ
 	
 	\noindent and
 	
 	\DEQSZ
 	\widehat{uv} = \widehat{u}\ast\widehat{v} \quad \text{and}\quad \widehat{u\ast v} = \widehat{u} \, \widehat{v}. \label{Fourier-prodoct}
 	\EEQSZ
 }
  {In \eqref{Fourier-derivative}, $x$ and $\xi$ in $\partial^\alpha_x$ and $\partial^\alpha_\xi$ indicate the $\alpha^{th}$ order derivative with respect to space variables in the Euclidean and Fourier spaces respectively, $\ast$ in \eqref{Fourier-prodoct} is the convolution operator and Equation \eqref{plancherel} is the Parseval-Plancherel identity. }
 For the detail of  these and other properties of the  Fourier transform we refer to \cite{bahouri2011fourier,hormander1983analysis,wolff2003lectures}.
 
 In fact, \eqref{plancherel} implies  that  energy of the system \eqref{mhdmain} in  Fourier space is equal to  energy of the system in Cartesian space.  
 To take advantage of \eqref{plancherel}  we give an equivalent formulation for  \eqref{mhdmain} in Fourier space. 
 %
 %
 %
  {This is done in two steps; first we eliminate the pressure term by applying the Leray projector given by \eqref{leray}. 
 }
 
 \DEQSZ
 {\CP}\cdot := Id -\nabla\Delta^{-1}{\rm div } \cdot. \label{leray} 
 \EEQSZ
 
 The application of   ${\CP}$  together with the fact that the fields $u$ and $v$ and the inhomogeneous terms $f_1$ and $f_2$ are divergence free reduces the system \eqref{mhdmain} to
 
 \DEQSZ
 \left\lbrace \begin{array}{l l}
 	\partial_tu - \nu\Delta u={\CP}((b\cdot \nabla)b)-{\CP}((u\cdot \nabla)u) + f_1,& \\
 	\partial_t b - \eta\Delta b={\CP}((b\cdot \nabla)u)-{\CP}((u\cdot \nabla)b) + f_2,&\\
 	u|_{t=0}=\,u_0 \quad b|_{t=0}=\,b_0 .&
 \end{array} \right.\label{lmhd} 
 \EEQSZ
 
 \noindent Next, we take the Fourier transform of \eqref{lmhd} to get
 
 \DEQSZ
 \left\lbrace \begin{array}{l l}
 	\widehat{u}_t + \nu |k|^2\widehat{u}=\reallywidehat{({\CP} ((b\cdot\nabla )b) )} -  \reallywidehat{({{\CP} ((u\cdot\nabla )u) })} +\widehat{f_1},&\\
 	\widehat{b}_t + \eta |k|^2\widehat{b}= \reallywidehat{({\CP} ((b\cdot\nabla )u) )} -  \reallywidehat{({{\CP} ((u\cdot\nabla )b) })} +\widehat{f_2} ,&\\
 	\widehat{u}|_{t=o}=\widehat{u}_0,\quad \widehat{b}|_{t=o}=\widehat{b}_0 .&
 \end{array} \right.   \label{flmhd}
 \EEQSZ
 
 \noindent Thus \eqref{flmhd} is 
 an equivalent formulation of \eqref{mhdmain} in Fourier space.

 \subsection{ {The Estimates }  }
 \label{fR3}
 
 This section is devoted to finding estimates in Fourier space for solutions of \eqref{mhdmain}. 
 For ease of calculations, we define an operator 
 
 \DEQSZ\label{operator: pi}
 \Pi_\xi:\mathbb{C}^3 \to\mathbb{C}^2_\xi  \quad \text{by} \quad \Pi_\xi(z) = z-(z\cdot \xi)\frac{\xi}{|\xi|^2},
 \EEQSZ
 
 \noindent where $ \mathbb{C}^3  $ the usual three dimensional complex space and 
  {\DEQS
 	\mathbb{C}^2_\xi :=\{z\in\mathbb{C}^3: \xi\cdot z = 0\}.
 	\EEQS
}
 
 Observe that for $\xi\in\mathbb{C}^3$ and $u$ divergence free, we have
 
 \DEQSZ
 \pi_\xi(\hat{u}) &=& \hat{u} \nonumber\\
 \reallywidehat{\CP((u\cdot \nabla)b)} &=& i\Pi_\xi\lk( \int_D\widehat{u}(\xi-\zeta)\zeta\widehat{b}(\zeta) {\rm d}\zeta   \rk) = i\Pi_\xi\lk( \int_D\zeta\widehat{u}(\xi-\zeta)\widehat{b}(\zeta) {\rm d}\zeta   \rk)  .\label{topik}
 \EEQSZ
 
 \noindent Now plugging \eqref{topik} in \eqref{flmhd} we get,
 
 \DEQSZ
 \left\lbrace \begin{array}{l l}
 	\widehat{u}_t =- \nu |\xi|^2\widehat{u}+ i\Pi_\xi(\int_D\zeta \widehat{b}(\xi-\zeta)\widehat{b}(\zeta)d\zeta) -i\Pi_\xi(\int_D\zeta \widehat{u}(\xi-\zeta)\widehat{u}(\zeta)d\zeta) +\widehat{f_1},&\\
 	\widehat{b}_t =- \eta |\xi|^2\widehat{b}+i\Pi_\xi(\int_D\zeta \widehat{b}(\xi-\zeta)\widehat{u}(\zeta)d\zeta) -i\Pi_\xi(\int_D\zeta \widehat{u}(\xi-\zeta)\widehat{b}(\zeta)d\zeta) +\widehat{f_2} ,& \\
 	\widehat{u}|_{t=o}=\widehat{u}_0,\quad \widehat{b}|_{t=o}=\widehat{b}_0. &
 \end{array}\right. \label{aaaa}
 \EEQSZ 
 
 
 
 

  {
 	\begin{remark}\label{assum: purposefully for 17}
 		Let $B_R(0)$  a ball in $L^2(D)$ of radius $R$. Let $(u_0,b_0)\in B_R(0)  $ and  $f_1,f_2\in L_{loc}^\infty([0,\infty);H^{-1}(D)\cap L^2(D))$. 
 		If   
 		an appropriate frame is chosen and the total pressure $\Pi$ is suitably normalized so that 
 		
 		\DEQSZ
 		\int_D u(x,t)\cdot f_1(x,t) + b(x,t)\cdot f_2(x,t) {\rm d} x, \label{assum: F(T)}
 		\EEQSZ
 		
 		\noindent is bounded, 
 		then  for any $T>0$  there is a non negative function $R(T)$ such that 
 		
 		\DEQSZ
 		\| u(\cdot,T) \|_{L^2}^2 + \| b(\cdot,T) \|_{L^2}^2+ \min(\nu,\eta)\int_0^T( \| \nabla u(\cdot,s) \|_{L^2}^2 + \| \nabla b(\cdot,s) \|_{L^2}^2){\rm d} s  
 		\leq R^2(T). \label{17}
 		\EEQSZ
 		
 		\noindent Furthermore, when  $f_1 \equiv f_2 \equiv 0 $, the bound $R(T) = R $ is a constant fully determined by the initial data $ (u_0,b_0) $.     In this case one could actually take $R$ to be the right hand side (RHS) of \eqref{lerrayenergy} and $B_R(0)$, a ball of radius $R$ and center $0,$ becomes an invariant \footnote{Set $A$ is said to be an invariant (future invariant) set with respect to a  function $\varphi$ or family of functions $\{\varphi(t): t\in\left[0,\infty \right)\}$, if
 			$$
 			\varphi(0) \in A \implies \varphi(t)\in A, \quad \forall t\geq 0.
 			$$}
 		set   for the weak solution.  
 		
 	\end{remark}
 }
 Assuming that the non-homogeneous terms $f_1$ and $f_2$ are appropriately chosen so that \eqref{17} holds. With no lose of generality, one may assume from 
 \eqref{plancherel} that 
 \DEQSZ
 \| \widehat{ u}(\cdot,t) \|_{L^2}^2+ \|  \widehat{b}(\cdot,t) \|_{L^2}^2\leq R^2(t).\label{eq: formainthm}
 \EEQSZ
 But the   problem   is since $u,\, b$  are only   distributional (weak)  solutions    their Fourier transforms are not well defined at particular points, say $(\xi,t),$ in Fourier-space-time.
 We  address the problem by taking a smooth cutoff of   $u$ and $b$ over a cube of finite length and making use of Paley-Wienner Theorem \cite[pp. 193]{hormander1983analysis}. 
 
 Let $k(\neq 0)\in\mathbb{R}^3,\, 0<\delta<\frac{|k|}{2\sqrt{3}}$.
 Define $\displaystyle {\chi}_k(\cdot)$ to be a smooth cutoff function of a cube $Q_k$ about $k$ of side length $2\delta$  
 such that 
 \DEQS
 \widehat{\chi}_k(\xi)=1,
 \EEQS 
 on a cube of the same center with side $\delta$ and 
 
 \DEQS
 {\rm supp} \widehat{\chi}=\{\xi\in\mathbb{R}^3: \frac{|k|}{2}\leq  |\xi|\leq\frac{3}{2}|k| \}.
 \EEQS
 
 \noindent Consider the following three smooth cutoff functions defined to suit our purpose; 
 
 \DEQSZ
 \left(\widehat{\chi}_k(D)u\right)(x,t)&:=& {\CF}^{-1}\left(\widehat{\chi}_k(\xi)\widehat{u}(\xi,t) \right) =(\chi_k*u)(x,t), \label{cutoff} \\
 %
 e_p(k,t)&:=&\left(\int_D| \widehat{\chi}_k(\xi)\widehat{u}(\xi,t)|^p+| \widehat{\chi}_k(\xi)\widehat{b}(\xi,t)|^p d\xi\right)^{\frac{1}{p}},  \label{ep} \\
 %
 h_p(k,t)&:=& \sup\limits_{0\leq s\leq t}\left(\int_D[|\widehat{\chi}_k(\xi)\widehat{f}_1(\xi,t) |^p + |\widehat{\chi}_k(\xi)\widehat{f}_1(\xi,t) |^p ]/{|\xi|^p}d\xi  \right)^{\frac{1}{p}}. \label{fp} 
 \EEQSZ

 \begin{remark} \label{Lipremark}
 	Since the Fourier transform of $\chi_k$ is compactly supported, by Paley-Wienner theorem, \cite[Theorem 7.3.1]{hormander1983analysis} we have $\chi_k\in H^m$ for all $m.$ Thus  $\chi_k$ can be considered as   a test function. 
 	
 \end{remark}
 We now have enough preparation to start working on  estimating our solution in Fourier space.  
 To establish necessary estimates, we first 
 need to establish estimates on $e_p,$ for $p=2$  
 followed by estimate for $e_p(k,t)$ for all $2\leq p\leq\infty $.  
 \begin{lemma}\label{lemma2.9}
 	Suppose that \eqref{17} holds
 	%
 	%
 	and there exists a non-decreasing function $R_1(t)$ such that
 	
 	\DEQSZ
 	(2\delta)^{3/2}\sqrt{2}R^2(t)+2h_2(k,t)<\frac{\min(\nu,\eta)}{6}R_1(t), \label{30}
 	\EEQSZ
 	
 	for all $t\in[0,\infty)$ and $\delta < \frac{|k|}{2\sqrt{3}}.$  
 	If $e_2(k,0)<\frac{R_1(0)}{|k|},$ then for any $t\in (0,\infty)$ we have
 	
 	\DEQSZ
 	e_2(k,t)\leq \frac{R_1(t)}{|k|} .\label{31}
 	\EEQSZ
 	
 \end{lemma}
 
 \begin{proof}[Proof of Lemma \ref{lemma2.9}]
 	By definition
 	
 	\DEQSZ
 	e_2^2(k,t)=\int_D \widehat{\chi}_k(\xi)\widehat{u}(\xi,t)\overline{\widehat{\chi}_k(\xi)\widehat{u}(\xi,t) }+  \widehat{\chi}_k(\xi)\widehat{b}(\xi,t)\overline{\widehat{\chi}_k(\xi)\widehat{b}(\xi,t)} {\rm d}\xi.\label{step1}
 	\EEQSZ
 	
 	Differentiating \eqref{step1} with respect to time and using equation \eqref{aaaa}, we get
 	
 	\DEQS
 	\lqq{\frac{d}{dt}e_2^2(k,t)} \\
 	&=&  \int_D \left[\widehat{\chi}_k(\xi) \frac{d}{dt}\lk(\widehat{u}(\xi,t)\rk)(\overline{\widehat{\chi}_k(\xi)\widehat{u}(\xi,t) }) + ( \widehat{\chi}_k(\xi)\widehat{u}(\xi,t)) \overline{\widehat{\chi}_k(\xi) }\frac{d}{dt}\lk(\overline{\widehat{u}(\xi,t) }\rk) \right. \nonumber \\
 	&&\left. + \widehat{\chi}_k(\xi)\frac{d}{dt}\lk(\widehat{b}(\xi,t)\rk)\overline{\widehat{\chi}_k(\xi)\widehat{b}(\xi,t)}+ \widehat{\chi}_k(\xi)\widehat{b}(\xi,t)\overline{ \widehat{\chi}_k(\xi)}\frac{d}{dt}\lk(\overline{\widehat{b}(\xi,t)}\rk)
 	\right]{\rm d}\xi \\
 	&=& \int_D\left[ \widehat{\chi}_k(\xi)\left(- \nu |\xi|^2\widehat{u}+ i\Pi_\xi(\int_D\zeta \widehat{b}(\xi-\zeta)\widehat{b}(\zeta){\rm d}\zeta)-\right.\right.\nonumber\\
 	&&\left.i\Pi_\xi(\int_D\zeta \widehat{u}(\xi-\zeta)\widehat{u}(\zeta){\rm d}\zeta) +\widehat{f_1}\right) (\overline{\widehat{\chi}_k(\xi)\widehat{u}(\xi,t) })\nonumber\\
 	&&+( \widehat{\chi}_k(\xi)\widehat{u}(\xi,t))\overline{\widehat{\chi}_k(\xi)\left(- \nu |\xi|^2\widehat{u}+ i\Pi_\xi(\int_D\zeta \widehat{b}(\xi-\zeta)\widehat{b}(\zeta){\rm d}\zeta)\right. } \nonumber\\
 	&& \overline{ \left. -i\Pi_\xi(\int_D\zeta \widehat{u}(\xi-\zeta)\widehat{u}(\zeta){\rm d}\zeta) +\widehat{f_1}\right)}\nonumber\\
 	&&+ \widehat{\chi}_k(\xi)\left( - \eta |\xi|^2\widehat{b}+i\Pi_\xi(\int_D\zeta \widehat{b}(\xi-\zeta)\widehat{u}(\zeta){\rm d}\zeta) \right. \nonumber \\
 	&& \left.-i\Pi_\xi(\int_D\zeta \widehat{u}(\xi-\zeta)\widehat{b}(\zeta){\rm d}\zeta) +\widehat{f_2}\right)\overline{\widehat{\chi}_k(\xi)\widehat{b}(\xi,t)}\nonumber\\
 	&&\left.+ \widehat{\chi}_k(\xi)\widehat{b}(\xi,t)\overline{\widehat{\chi}_k(\xi) - \eta |\xi|^2\widehat{b}+i\Pi_\xi(\int_D\zeta \widehat{b}(\xi-\zeta)\widehat{u}(\zeta){\rm d}\zeta)  {\rm d}\xi.} \right.  \nonumber \\
 	&&\left.-\overline{i\Pi_\xi(\int_D\zeta \widehat{u}(\xi-\zeta)\widehat{b}(\zeta){\rm d}\zeta) +\widehat{f_2} }\right].
 	\EEQS
 	
 	Applying elementary properties of complex numbers, it follows that
 	
 	\DEQSZ
 	\lqq{\frac{1}{2}\frac{d}{dt}e_2^2(k,t)}\label{32}\\
 	&=& -\nu\int_D |\xi|^2|\widehat{\chi}_k(\xi)\widehat{u}(\xi,t)|^2 {\rm d}\xi - \eta\int_D |\xi|^2|\widehat{\chi}_k(\xi)\widehat{b}(\xi,t)|^2{\rm d}\xi\nonumber\\
 	&&+\int_D\Re\left(i\Pi_\xi\left( \int_D\widehat{u}(\xi-\zeta)\cdot\zeta\widehat{u}(\zeta){\rm d}\zeta\right)\overline{\widehat{\chi}_k(\xi)\widehat{u}(\xi,t)}\right){\rm d}\xi\nonumber\\
 	&&+\int_D\Re\left(i\Pi_\xi\left( \int_D\widehat{b}(\xi-\zeta)\cdot\zeta\widehat{u}(\zeta){\rm d}\zeta\right)\overline{\widehat{\chi}_k(\xi)\widehat{b}(\xi,t)}\right){\rm d}\xi\nonumber \\
 	&& -\int_D\Re\left(i\Pi_\xi\left( \int_D\widehat{b}(\xi-\zeta)\cdot\zeta\widehat{b}(\zeta){\rm d}\zeta\right)\overline{\widehat{\chi}_k(\xi)\widehat{u}(\xi,t)}\right){\rm d}\xi \nonumber\\
 	&&-\int_D\Re\left(i\Pi_\xi\left( \int_D\widehat{u}(\xi-\zeta)\cdot\zeta\widehat{b}(\zeta){\rm d}\zeta\right)\overline{\widehat{\chi}_k(\xi)\widehat{b}(\xi,t)}\right){\rm d}\xi \nonumber\\
 	&&+\Re\int_D  \widehat{\chi}_k(\xi)\widehat{f}_1(\xi, t)\overline{\widehat{\chi}_k(\xi)\widehat{u}(\xi,t)}{\rm d}\xi +\Re\int_D \widehat{\chi}_k(\xi)\widehat{f}_2(\xi, t)\overline{\widehat{\chi}_k(\xi)\widehat{b}(\xi,t)}{\rm d}\xi\nonumber\\
 	&:=& I_1 + I_2+I_3+I_4+I_5+I_6+I_7+I_8. \nonumber
 	\EEQSZ
 	
 	For ease of calculations, we now deal  with the terms on RHS of \eqref{32} separately.

 	\DEQSZ
 	I_1 +I_2 &= -\nu\int_D |\xi|^2|\widehat{\chi}_k(\xi)\widehat{u}(\xi,t)|^2 {\rm d}\xi - \eta\int_D |\xi|^2|\widehat{\chi}_k(\xi)\widehat{b}(\xi,t)|^2 {\rm d}\xi \nonumber\\ 
 	&\leq -\min(\nu,\eta)\frac{|k|^2}{4}\int_D\left(|\widehat{\chi}_k(\xi)\widehat{u}(\xi,t)|^2 + |\widehat{\chi}_k(\xi)\widehat{b}(\xi,t)|^2\right){\rm d}\xi \label{I12est}.
 	\EEQSZ
 	In \eqref{I12est} we used the fact $ \xi\in{\rm supp }\widehat{\chi}_k $; that is $ \frac{|k|}{2}\leq|\xi|\leq\frac{3}{2}|k|. $
 	
 	\DEQS
 	I_3 = -\Im\int_D\left(\Pi_\xi\left( \int_D\widehat{u}(\xi-\zeta)\cdot\zeta\widehat{u}(\zeta){\rm d}\zeta\right)\overline{\widehat{\chi}_k(\xi)\widehat{u}(\xi,t)}\right){\rm d}\xi, 
 	\EEQS
 	
 	which implies
 	
 	\DEQSZ
 	|I_3| &\leq & \left|\int_D\left(\Pi_\xi\left( \int_D\widehat{u}(\xi-\zeta)\cdot\zeta\widehat{u}(\zeta){\rm d}\zeta\right)\overline{\widehat{\chi}_k(\xi)\widehat{u}(\xi,t)}\right){\rm d}\xi \right| \label{step2}\\
 	&\leq&  \| \widehat{\chi}_k\widehat{u}(\cdot,t) \|_{L^2}  \| \widehat{\chi}_k \Pi_\xi\left( \xi\cdot\int_D\widehat{u}(\xi-\zeta)\widehat{u}(\zeta){\rm d}\zeta \right)  \|_{L^2} \nonumber\\
 	&\leq &  \| \widehat{\chi}_k\widehat{u}(\cdot,t) \|_{L^2}  \| \xi\widehat{\chi}_k \|_{L^2} \|\widehat{u}(\cdot,t)  \|_{L^2}^2. \nonumber
 	\EEQSZ
 	
 	The estimate in \eqref{step2} is due to the fact that     $u\, \text{and}\,b$ are divergence free and elementary properties of complex numbers.  {H\"older's}  and Young's  inequalities are also used.
 	
 	We know from construction of $\chi_k$ and H\"older's inequality that 
 	
 	\DEQSZ
 	\| \xi\widehat{\chi}_k \|_{L^2} &\leq& \| \xi \|_{L^4} \|\chi_k \|_{L^4} \label{step3} \\
 	&=&  \left(\int_{Q_k}|\xi|^4{\rm d}\xi\right)^{\frac{1}{4}}\left(\int_{Q_k}|\chi_k|^4{\rm d}\xi\right)^{\frac{1}{4}} \nonumber\\
 	&\leq& \frac{3}{2}|k|(2\delta)^\frac{3}{4} (2\delta)^\frac{3}{4} = \frac{3}{2}|k|(2\delta)^\frac{3}{2}.\nonumber 
 	\EEQSZ
 	
 	Thus, combining \eqref{step2} and \eqref{step3} we get,
 	
 	\DEQSZ
 	|I_3| \leq\frac{3}{2}|k|(2\delta)^\frac{3}{2} \| \widehat{\chi}_k\widehat{u}(\cdot,t) \|_{L^2} \|\widehat{u}(\cdot,t)  \|_{L^2}^2.\label{I3est}
 	\EEQSZ
 	
 	Proceeding similarly with $I_4,\, I_5$ and $I_6$ we get
 	
 	\DEQSZ
 	|I_4| &\leq& \frac{3}{2}|k|(2\delta)^\frac{3}{2} \| \widehat{\chi}_k\widehat{b}(\cdot,t) \|_{L^2} \|\widehat{u}(\cdot,t)  \|_{L^2} \|\widehat{b}(\cdot,t)  \|_{L^2} \label{I4est},\\
 	|I_5| &\leq & \frac{3}{2}|k|(2\delta)^\frac{3}{2} \| \widehat{\chi}_k\widehat{u}(\cdot,t) \|_{L^2} \|\widehat{b}(\cdot,t)  \|_{L^2}^2 \label{I5est},\\
 	|I_6| &\leq & \frac{3}{2}|k|(2\delta)^\frac{3}{2} \| \widehat{\chi}_k\widehat{b}(\cdot,t) \|_{L^2} \|\widehat{u}(\cdot,t)  \|_{L^2} \|\widehat{b}(\cdot,t)  \|_{L^2} .\label{I6est}
 	\EEQSZ
 	
 	Thanks to H\"older's inequality, the integral   $I_7$ is estimated as follows;
 	
 	\DEQSZ
 	|I_7 | &=&  \left| \Re\int_D  \widehat{\chi}_k(\xi)\widehat{f}_1(\xi, t)\overline{\widehat{\chi}_k(\xi)\widehat{u}(\xi,t)}d\xi \right|  \label{I7est}\\
 	&\leq & \left|\int_D \overline{\widehat{\chi}_k(\xi)\widehat{u}(\xi,t)} \widehat{\chi}_k(\xi)|\xi|\widehat{f}_1(\xi, t)/{|\xi|}d\xi\right|\nonumber\\
 	&\leq & \| |\xi|\widehat{\chi}_k\widehat{u}(\cdot,t) \|_{L^2} \| \widehat{\chi}_k\widehat{f}_1(\cdot,t){|\xi|^{-1}} \|_{L^2} \nonumber \\
 	&\leq & { \frac{3}{2}|k| \| \widehat{\chi}_k\widehat{u}(\cdot,t) \|_{L^2}  \| \widehat{\chi}_k\widehat{f}_1(\cdot,t){|\xi|^{-1}} \|_{L^2}}.\nonumber
 	\EEQSZ 
 	
 	Similarly, we have
 	
 	\DEQSZ
 	|I_8 | &\leq & { \frac{3}{2}|k| \| \widehat{\chi}_k\widehat{b}(\cdot,t) \|_{L^2}  \| \widehat{\chi}_k\widehat{f}_2(\cdot,t){|\xi|^{-1}} \|_{L^2}} \label{I8est}.
 	\EEQSZ
 	
 	Now combining the estimates \eqref{I12est}-\eqref{I8est} we obtain
 	
 	\DEQSZ 
 	\lqq{\frac{1}{2}\frac{d}{dt}e_2^2(k,t)} \label{step4}\\
 	&=&-\min(\nu,\eta)\frac{|k|^2}{4}\int_D\left(|\widehat{\chi}_k(\xi)\widehat{u}(\xi,t)|^2 + |\widehat{\chi}_k(\xi)\widehat{b}(\xi,t)|^2\right)d\xi \nonumber \\
 	&&+ \frac{3}{2}|k|(2\delta)^\frac{3}{2} \| \widehat{\chi}_k\widehat{u}(\cdot,t) \|_{L^2} \left( \|\widehat{u}(\cdot,t)  \|_{L^2}^2 +  \|\widehat{b}(\cdot,t)  \|_{L^2}^2 \right) \nonumber\\
 	&&+ 3|k|(2\delta)^\frac{3}{2} \| \widehat{\chi}_k\widehat{b}(\cdot,t) \|_{L^2} \|\widehat{u}(\cdot,t)  \|_{L^2} \|\widehat{b}(\cdot,t)  \|_{L^2} \nonumber\\
 	&&+  \frac{3}{2}|k|\left[ \| \widehat{\chi}_k\widehat{u}(\cdot,t) \|_{L^2}  \| \widehat{\chi}_k\widehat{f}_1(\cdot,t){|\xi|^{-1}} \|_{L^2} +  \| \widehat{\chi}_k\widehat{b}(\cdot,t) \|_{L^2}  \| \widehat{\chi}_k\widehat{f}_2(\cdot,t){|\xi|^{-1}} \|_{L^2} \right]\nonumber \\
 	&\leq& -\min(\nu,\eta)\frac{|k|^2}{4}e_2^2(k,t) + \frac{3}{2}|k|(2\delta)^{3/2}\left( \|\widehat{\chi}_k\widehat{u} \| +  \|\widehat{\chi}_k\widehat{b} \|\right)\left( \|\widehat{u} \|^2+ \|\widehat{b} \|^2 \right) \nonumber\\
 	&& + \frac{3}{2}|k|\left( \|\widehat{\chi}_k\widehat{u} \| +  \|\widehat{\chi}_k\widehat{b} \|\right)\left(  \|\widehat{\chi}_k\widehat{f}_1/|\xi| \| +  \|\widehat{\chi}_k\widehat{f}_2/|\xi| \| \right) \nonumber\\
 	&\leq & -\min(\nu,\eta)\frac{|k|^2}{4}e_2^2(k,t) + \frac{3}{2}|k|(2\delta)^{3/2}\sqrt{2}e_2(k,t)\left( \|\widehat{u} \|^2+  \|\widehat{b} \|^2 \right)\nonumber \\
 	&& + \frac{3}{2}|k|\sqrt{2}e_2(k,t)\left(  \|\widehat{\chi}_k\widehat{f}_1/|\xi| \| +  \|\widehat{\chi}_k\widehat{f}_2/|\xi| \| \right) \nonumber\\
 	&\leq & -\min(\nu,\eta)\frac{|k|^2}{4}e_2^2(k,t) +\frac{3}{2}|k|e_2(k,t)\left((2\delta)^{3/2}\sqrt{2}R^2(t) + 2h_2(k,t) \right) \nonumber.
 	\EEQSZ
 	
 	Here we used Serine's inequality  \cite[Lemma 1]{serrin1964local} to estimate  upper bounds for $  \|\widehat{\chi}_k\widehat{u} \| +  \|\widehat{\chi}_k\widehat{b} \| $ and  
 	$  \|\widehat{\chi}_k\widehat{f}_1/|\xi| \| +  \|\widehat{\chi}_k\widehat{f}_2/|\xi| \| $ respectively as;
 	
 	\DEQS
 	&& \|\widehat{\chi}_k\widehat{u} \| +  \|\widehat{\chi}_k\widehat{b} \|\leq \sqrt{2}e_2(k,t) ,\\
 	\\
 	&& \|\widehat{\chi}_k\widehat{f}_1/|\xi| \| +  \|\widehat{\chi}_k\widehat{f}_2/|\xi| \| \leq \sqrt{2}h_2.
 	\EEQS 
 	
 	Now define the set $B_{R_1}$ by, 
 	
 	\DEQSZ
 	B_{R_1}=\left\lbrace e:e\leq R_1/{|k|} \right\rbrace = \left\lbrace e(k,t): e(k,t)\leq R_1(t)/|k| \right\rbrace . \label{attract}
 	\EEQSZ		
 	
 	When $e(k,t)=e_2(k,t)=\frac{R_1(t)}{|k|}$ in \eqref{attract}, we get 
 	
 	\DEQS
 	&\frac{1}{2}\frac{d}{dt}e_2^2(k,t) < \frac{-\min(\nu,\eta)}{4}R_1^2(t)+\frac{3}{2}R_1(t)\frac{\min(\nu,\eta)}{6}R_1(t) \leq 0. 
 	\EEQS
 	
 	Then by chain rule and
 	from the fact that  $e_2(k,t)\geq0$, we conclude that 
 	
 	\DEQSZ
 	\frac{d}{dt}e_2(k,t) <0. \label{tocon}
 	\EEQSZ
 	
 	Indeed, \eqref{tocon} implies that 
 	$B_{R_1}$ is an attracting    set for $e_2(k,t)$.
 	Therefore, if $e_2(k,0)<\frac{R_1(0)}{|k|},$ then $e_2(k,t)<\frac{R_1(t)}{|k|}$ for all $t\in (0,\infty).$ 
 \end{proof}
 
 \begin{lemma}\label{lemma2.10}
 	Suppose that for a given $k\in\mathbb{R}^3$ and $2\leq p<\infty$ there is a nondecreasing function $R_1(t)$ that satisfies the condition
 	
 	\DEQS
 	2^\frac{1}{p}(2\delta)^{3/p}R^2(t) +  2h_p(k,t)<\frac{\min(\nu,\eta)}{6}R_1(t),
 	\EEQS
 	
 	for  $0<\delta<|k|/{2\sqrt{3}}.$
 	
 	If a solution to \eqref{mhdmain} initially satisfies 
 	
 	\DEQS
 	e_p(k,0)<R_1(0)/|k|,
 	\EEQS
 	
 	then for all $0<t<\infty$, 
 	
 	\DEQS 
 	e_p(k,t)<\frac{R_1(t)}{|k|}.
 	\EEQS 
 	
 \end{lemma}
 
 \begin{proof}[Proof of Lemma \ref{lemma2.10}]
 	The proof follows same procedure as the proof of Lemma \ref{lemma2.9}.
 	We begin by  taking the time derivative of $e_p^p(k,t)$.		
 	
 	\DEQSZ
 	\lqq{\frac{d}{dt}e_p^p(k,t)}\label{depp}\\
 	&=&\partial_t\int\left|\widehat{\chi}_k(\xi)\widehat{u}(\xi,t)\right|^p+\left|\widehat{\chi}_k(\xi)\widehat{b}(\xi,t)\right|^p {\rm d}\xi \nonumber\\
 	&=&\Re\left\lbrace \int\left(p|\widehat{\chi}_k(\xi) \widehat{u}(\xi,t)|^{p-2}\left((\widehat{\chi}_k(\xi)\partial_t\widehat{u}(\xi,t))(\overline{\widehat{\chi}_k(\xi) \widehat{u}(\xi,t)}) \right)\right.\right.\nonumber\\
 	&&\left.\left.+p|\widehat{\chi}_k(\xi) \widehat{b}(\xi,t)|^{p-2}\left((\widehat{\chi}_k(\xi)\partial_t\widehat{b}(\xi,t))(\overline{\widehat{\chi}_k(\xi) \widehat{b}(\xi,t)}) \right)\right){\rm d}\xi  \right\rbrace\nonumber\\
 	&=&-\nu\int p|\xi|^2|\widehat{\chi}_k(\xi) \widehat{u}(\xi,t)|^p{\rm d}\xi-\eta\int p|\xi|^2|\widehat{\chi}_k(\xi) \widehat{b}(\xi,t)|^p{\rm d}\xi \nonumber\\
 	&&+\Re\int ip|\widehat{\chi}_k(\xi) \widehat{u}(\xi,t)|^{p-2}\overline{\widehat{\chi}_k(\xi) \widehat{u}(\xi,t) }\widehat{\chi}_k(\xi)\Pi_\xi\int\widehat{u}(\xi-\zeta) \zeta\widehat{u}(\zeta){\rm d}\zeta {\rm d}\xi \nonumber\\
 	&&+\Re\int ip|\widehat{\chi}_k(\xi) \widehat{u}(\xi,t)|^{p-2}\overline{\widehat{\chi}_k(\xi) \widehat{u}(\xi,t) }\widehat{\chi}_k(\xi)\Pi_\xi\int\widehat{b}(\xi-\zeta) \zeta\widehat{b}(\zeta){\rm d}\zeta {\rm d}\xi \nonumber\\
 	&&+\Re\int ip|\widehat{\chi}_k(\xi) \widehat{b}(\xi,t)|^{p-2}\overline{\widehat{\chi}_k(\xi) \widehat{b}(\xi,t) }\widehat{\chi}_k(\xi)\Pi_\xi\int\widehat{b}(\xi-\zeta) \zeta\widehat{u}(\zeta){\rm d}\zeta {\rm d}\xi \nonumber\\
 	&&+\Re\int ip|\widehat{\chi}_k(\xi) \widehat{b}(\xi,t)|^{p-2}\overline{\widehat{\chi}_k(\xi) \widehat{b}(\xi,t) }\widehat{\chi}_k(\xi)\Pi_\xi\int\widehat{u}(\xi-\zeta) \zeta\widehat{b}(\zeta){\rm d}\zeta {\rm d}\xi  \nonumber\\
 	&&+\Re\int p|\widehat{\chi}_k(\xi) \widehat{u}(\xi,t)|^{p-2}\overline{\widehat{\chi}_k(\xi) \widehat{u}(\xi,t) }\widehat{\chi}_k(\xi)\widehat{f}_1(\xi,t){\rm d}\xi  \nonumber\\
 	&&+\Re\int p|\widehat{\chi}_k(\xi) \widehat{b}(\xi,t)|^{p-2}\overline{\widehat{\chi}_k(\xi) \widehat{b}(\xi,t) }\widehat{\chi}_k(\xi)\widehat{f}_2(\xi,t){\rm d}\xi \nonumber \\
 	&=:& I_1+I_2+I_3+I_4+I_5+I_6+I_7+I_8 \nonumber.
 	\EEQSZ
 	
 	In  the derivation of \eqref{depp} we have used the following fact;
 	
 	\DEQS
 	\frac{{\rm d}}{{\rm d} t}|\widehat{\chi}_k(\xi)\widehat{b}(\xi,t) | &=&	\frac{d}{dt}\sqrt{\widehat{\chi}_k(\xi)\widehat{b}(\xi,t)\overline{ \widehat{\chi}_k(\xi)\widehat{b}(\xi,t)}} \\
 	&=& \frac{1}{2}|\widehat{\chi}_k(\xi)\widehat{b}(\xi,t)|^{-1}\left((\widehat{\chi}_k(\xi)\partial_t\widehat{b}(\xi,t))(\overline{\widehat{\chi}_k(\xi) \widehat{b}(\xi,t)}) \right.\\  
 	&&\left.+ \overline{(\widehat{\chi}_k(\xi)\partial_t\widehat{b}(\xi,t))(\overline{\widehat{\chi}_k(\xi) \widehat{b}(\xi,t)}) } \right).
 	\EEQS
 	
 	We now estimate the integrals at the RHS of \eqref{depp}.
 	
 	\DEQSZ
 	I_1+I_2 &=& -\nu\int p|\xi|^2|\widehat{\chi}_k(\xi) \widehat{u}(\xi,t)|^p{\rm d}\xi-\eta\int p|\xi|^2|\widehat{\chi}_k(\xi) \widehat{b}(\xi,t)|^p{\rm d}\xi \nonumber \\
 	&\leq& \frac{-\nu p|k|^2}{4}\int|\widehat{\chi}_k(\xi) \widehat{u}(\xi,t)|^p {\rm d}\xi +\frac{-\eta p|k|^2}{4}\int|\widehat{\chi}_k(\xi) \widehat{b}(\xi,t)|^p {\rm d}\xi\nonumber \\
 	& \leq & -\frac{\min(\nu,\eta)}{4}p|k|^2e_p^p(k,t). \label{I12est2} 
 	\EEQSZ
 	
 	Here we   used the fact that for $\xi\in{\rm supp }\,\widehat{\chi}$, $ \frac{|k|}{2}\leq |\xi|\leq\frac{3}{2}|k|.$
 	Finally, thanks to H\"{o}lder's and Young's inequalities, we have 
 	
 	\DEQSZ
 	|I_3 |&=& \left|\Im\int p|\widehat{\chi}_k(\xi) \widehat{u}(\xi,t)|^{p-2}\overline{\widehat{\chi}_k(\xi) \widehat{u}(\xi,t) }\widehat{\chi}_k(\xi)\Pi_\xi\xi \int\widehat{u}(\xi-\zeta)\widehat{u}(\zeta){\rm d}\zeta {\rm d}\xi \right|\label{I3est2}\\
 	&\leq& \left|\int p|\widehat{\chi}_k(\xi) \widehat{u}(\xi,t)|^{p-2}\overline{\widehat{\chi}_k(\xi) \widehat{u}(\xi,t) }\widehat{\chi}_k(\xi)\Pi_\xi\xi \int\widehat{u}(\xi-\zeta)\widehat{u}(\zeta){\rm d}\zeta {\rm d}\xi \right| \nonumber\\
 	&\leq& p\left(\int\left( |\widehat{\chi}_k(\xi) \widehat{u}(\xi,t)|^{p-1}\right)^{\frac{p}{p-1}}{\rm d}\xi \right)^{\frac{p-1}{p}}\left(\int| \widehat{\chi}_k(\xi)\Pi_\xi\xi \int\widehat{u}(\xi-\zeta)\widehat{u}(\zeta)d\zeta|^p{\rm d}\xi \right)^{\frac{1}{p}}\nonumber \\
 	&\leq& p\left(\int\left( |\widehat{\chi}_k(\xi) \widehat{u}(\xi,t)|^{p-1}\right)^{\frac{p}{p-1}}{\rm d}\xi \right)^{\frac{p-1}{p}}\left( \int|\xi\widehat{\chi
 	}_k(\xi)|^p{\rm d}\xi\right)^{\frac{1}{p}}  \| \int\widehat{u}(\xi-\zeta)\widehat{u}(\zeta){\rm d}\zeta \|_{L^\infty} \nonumber\\
 	&\leq &p\left(\int\left( |\widehat{\chi}_k(\xi) \widehat{u}(\xi,t)|^{p-1}\right)^{\frac{p}{p-1}}{\rm d}\xi \right)^{\frac{p-1}{p}}\left( \int |\xi|^p|\widehat{\chi}_k(\xi)|^p{\rm d}\xi\right)^{\frac{1}{p}}\left(\int|\widehat{u}(\xi,t)|^p{\rm d}\xi\right)^{\frac{1}{p}}\nonumber\\
 	&&\times\left(\int|\widehat{u}(\xi,t)|^{\frac{p}{p-1}}{\rm d}\xi\right)^{\frac{p-1}{p}} \nonumber.
 	\EEQSZ
 	
 	Following a similar approach yields,
 	
 	\DEQSZ
 	|I_4| &\leq & p\left(\int\left( |\widehat{\chi}_k(\xi) \widehat{u}(\xi,t)|^{p-1}\right)^{\frac{p}{p-1}}{\rm d}\xi \right)^{\frac{p-1}{p}}\left( \int |\xi|^p|\widehat{\chi}_k(\xi)|^p{\rm d}\xi\right)^{1/p}\label{I4est2} \\
 	&&\times\left(\int|\widehat{b}(\xi,t)|^p{\rm d}\xi\right)^{1/p}\left(\int|\widehat{b}(\xi,t)|^{\frac{p}{p-1}}{\rm d}\xi\right)^{\frac{p-1}{p}} \nonumber,\\
 	|I_5|&\leq & p\left(\int\left( |\widehat{\chi}_k(\xi) \widehat{b}(\xi,t)|^{p-1}\right)^{\frac{p}{p-1}}{\rm d}\xi \right)^{\frac{p-1}{p}}\left( \int |\xi|^p|\widehat{\chi}_k(\xi)|^p{\rm d}\xi\right)^{1/p}\label{I5est2} \\
 	&&\times\left(\int|\widehat{b}(\xi,t)|^p{\rm d}\xi\right)^{1/p}\left(\int|\widehat{u}(\xi,t)|^{\frac{p}{p-1}}{\rm d}\xi\right)^{\frac{p-1}{p}} \nonumber,\\
 	|I_6|&\leq & p\left(\int\left( |\widehat{\chi}_k(\xi) \widehat{b}(\xi,t)|^{p-1}\right)^{\frac{p}{p-1}}{\rm d}\xi \right)^{\frac{p-1}{p}}\left( \int |\xi|^p|\widehat{\chi}_k(\xi)|^p{\rm d}\xi\right)^{1/p}\label{I6est2}\\
 	&&\times\left(\int|\widehat{u}(\xi,t)|^p{\rm d}\xi\right)^{1/p}\left(\int|\widehat{b}(\xi,t)|^{\frac{p}{p-1}}{\rm d}\xi\right)^{\frac{p-1}{p}}  \nonumber .
 	\EEQSZ
 	
 	We now remain to estimate $I_7$ and $I_8$. 
 	
 	\DEQSZ
 	|I_7| &=& \left| \Re\int p|\widehat{\chi}_k(\xi) \widehat{u}(\xi,t)|^{p-2}\overline{\widehat{\chi}_k(\xi) \widehat{u}(\xi,t) }\widehat{\chi}_k(\xi)\widehat{f}_1(\xi,t){\rm d}\xi \right|\label{I7est2}\\
 	&\leq&\left| \int p|\widehat{\chi}_k(\xi) \widehat{u}(\xi,t)|^{p-2}\overline{\widehat{\chi}_k(\xi) \widehat{u}(\xi,t) }\widehat{\chi}_k(\xi)\widehat{f}_1(\xi,t){\rm d}\xi \right|\nonumber\\
 	&\leq& p\left(\int|\widehat{\chi}_k(\xi) \widehat{u}(\xi,t)|^p{\rm d}\xi\right)^{\frac{p-1}{p}} \left(\int\frac{|\xi|^p|\widehat{\chi}_k(\xi)\widehat{f}_1(\xi,t)|^p}{|\xi|^p}{\rm d}\xi\right)^p\nonumber\\
 	&\leq& \frac{3p}{2}|k|\left(\int|\widehat{\chi}_k(\xi) \widehat{u}(\xi,t)|^p{\rm d}\xi\right)^{\frac{p-1}{p}} \left( \int\frac{ |\widehat{\chi}_k(\xi)\widehat{f}_1(\xi,t)|^p}{|\xi|^p }{\rm d}\xi\right)^{1/p} \nonumber.  
 	\EEQSZ
 	
 	And a similar approach yields, 
 	
 	\DEQSZ
 	|I_8|\leq \frac{3p}{2}|k|\left(\int|\widehat{\chi}_k(\xi) \widehat{b}(\xi,t)|^p{\rm d}\xi\right)^{\frac{p-1}{p}} \left( \int\frac{ |\widehat{\chi}_k(\xi)\widehat{f}_2(\xi,t)|^p}{|\xi|^p }{\rm d}\xi\right)^{1/p} \label{I8est2} .
 	\EEQSZ
 	
 	Now plugging the estimates \eqref{I12est2}-\eqref{I8est2} in \eqref{depp} and rearranging the terms we get,
 	
 	\DEQS
 	\frac{d}{dt}e_p^p(k,t) &\leq&-\frac{\min(\nu,\eta)}{4}p|k|^2e_p^p(k,t) \\
 	&&+\left( \int |\xi|^p|\widehat{\chi}_k(\xi)|^p{\rm d}\xi\right)^{1/p}\left[\left(\int |\widehat{\chi}_k(\xi)\widehat{u}(\xi,t)|^p{\rm d}\xi \right)^\frac{p-1}{p} \| \widehat{u} \|_{L^2}^2   \right. \\
 	&&+\left(\int |\widehat{\chi}_k(\xi)\widehat{u}(\xi,t)|^p{\rm d}\xi \right)^\frac{p-1}{p} \| \widehat{b} \|_{L^2}^2  
 	+\left. \left(\int |\widehat{\chi}_k(\xi)\widehat{b}(\xi,t)|^p{\rm d}\xi \right)^\frac{p-1}{p} \| \widehat{u} \|_{L^2}  \| \widehat{b} \|_{L^2} \right. \\
 	&& \left. +\left(\int |\widehat{\chi}_k(\xi)\widehat{b}(\xi,t)|^p{\rm d}\xi \right)^\frac{p-1}{p} \| \widehat{u} \|_{L^2}  \| \widehat{b} \|_{L^2}\right] \\
 	&&+ \frac{3p}{2}|k|\left(\int|\widehat{\chi}_k(\xi) \widehat{u}(\xi,t)|^p{\rm d}\xi\right)^{\frac{p-1}{p}} \left( \int\frac{ |\widehat{\chi}_k(\xi)\widehat{f}_1(\xi,t)|^p}{|\xi|^p }\right)^{1/p}\\
 	&&+\frac{3p}{2}|k|\left(\int|\widehat{\chi}_k(\xi) \widehat{b}(\xi,t)|^p{\rm d}\xi\right)^{\frac{p-1}{p}} \left( \int\frac{ |\widehat{\chi}_k(\xi)\widehat{f}_2(\xi,t)|^p}{|\xi|^p }{\rm d}\xi\right)^{1/p} .
 	\EEQS
 	
 	We know from the property of $\widehat{\chi}_k$ that $ \left( \int |\xi|^p|\widehat{\chi}_k(\xi)|^p{\rm d}\xi\right)^{1/p}$ is bounded from above as 
 	
 	\DEQSZ
 	\left( \int |\xi|^p|\widehat{\chi}_k(\xi)|^p{\rm d}\xi\right)^{1/p} \leq \frac{3|k|}{2}(2\delta)^{3/p}. \label{one} 
 	\EEQSZ 
 	
 	Furthermore, we have
 	
 	\DEQSZ
 	\lqq{\frac{3|k|p}{2}\left(\int|\widehat{\chi}_k(\xi) \widehat{u}(\xi,t)|^p{\rm d}\xi\right)^{\frac{p-1}{p}} \left( \int\frac{ |\widehat{\chi}_k(\xi)\widehat{f}_1(\xi,t)|^p}{|\xi|^p }\right)^{1/p} }
 	\label{two}\\ 
 	& +&\frac{3|k|p}{2}\left(\int|\widehat{\chi}_k(\xi) \widehat{b}(\xi,t)|^p{\rm d}\xi\right)^{\frac{p-1}{p}} \left( \int\frac{ |\widehat{\chi}_k(\xi)\widehat{f}_2(\xi,t)|^p}{|\xi|^p }{\rm d}\xi\right)^{1/p}\nonumber\\
 	&\leq & \frac{3|k|p}{2}\left(\left(\int|\widehat{\chi}_k(\xi) \widehat{u}(\xi,t)|^p{\rm d}\xi\right)^{\frac{p-1}{p}} + \left(\int|\widehat{\chi}_k(\xi) \widehat{b}(\xi,t)|^p{\rm d}\xi\right)^{\frac{p-1}{p}} \right) \nonumber \\
 	&&\times\left(\left( \int\frac{ |\widehat{\chi}_k(\xi)\widehat{f}_1(\xi,t)|^p}{|\xi|^p }\right)^{1/p} + \left( \int\frac{ |\widehat{\chi}_k(\xi)\widehat{f}_2(\xi,t)|^p}{|\xi|^p }\right)^{1/p}  {\rm d}\xi\right) \nonumber\\
 	&\leq & \frac{3|k|p}{2}2^\frac{1}{p}\left(\int|\widehat{\chi}_k(\xi) \widehat{u}(\xi,t)|^p{\rm d}\xi +\int|\widehat{\chi}_k(\xi) \widehat{b}(\xi,t)|^p{\rm d}\xi  \right)^{\frac{p-1}{p}}\nonumber \\
 	&& \times2^\frac{p-1}{p}\left(\int\frac{ |\widehat{\chi}_k(\xi)\widehat{f}_1(\xi,t)|^p}{|\xi|^p } +  \int\frac{ |\widehat{\chi}_k(\xi)\widehat{f}_2(\xi,t)|^p}{|\xi|^p }  {\rm d}\xi\right)^{1/p} \nonumber\\
 	&\leq & 2\frac{3|k|p}{2}e_p^{p-1}(k,t)h_p(k,t), \nonumber
 	\EEQSZ
 	
 	and 
 	
 	\DEQSZ
 	\lqq{\left[\left(\int |\widehat{\chi}_k(\xi)\widehat{u}(\xi,t)|^p{\rm d}\xi \right)^\frac{p-1}{p} \| \widehat{u} \|_{L^2}^2 +\left(\int |\widehat{\chi}_k(\xi)\widehat{u}(\xi,t)|^p{\rm d}\xi \right)^\frac{p-1}{p} \| \widehat{b} \|_{L^2}^2 \right.}  \label{three}\\
 	& +&\left. \left(\int |\widehat{\chi}_k(\xi)\widehat{b}(\xi,t)|^p{\rm d}\xi \right)^\frac{p-1}{p} \| \widehat{u} \|_{L^2}  \| \widehat{b} \|_{L^2} +\left(\int |\widehat{\chi}_k(\xi)\widehat{b}(\xi,t)|^p{\rm d}\xi \right)^\frac{p-1}{p} \| \widehat{u} \|_{L^2}  \| \widehat{b} \|_{L^2}\right]  \nonumber\\
 	&=& \left(\int |\widehat{\chi}_k(\xi)\widehat{u}(\xi,t)|^p{\rm d}\xi \right)^\frac{p-1}{p}\left( \|\widehat{u} \|^2+ \|\widehat{b} \|^2 \right) + 2\left(\int |\widehat{\chi}_k(\xi)\widehat{b}(\xi,t)|^p{\rm d}\xi \right)^\frac{p-1}{p} \|\widehat{u} \|  \|\widehat{b} \| \nonumber\\
 	&\leq& \left(\left(\int |\widehat{\chi}_k(\xi)\widehat{u}(\xi,t)|^p{\rm d}\xi \right)^\frac{p-1}{p} + \left(\int |\widehat{\chi}_k(\xi)\widehat{b}(\xi,t)|^p{\rm d}\xi \right)^\frac{p-1}{p}\right)\left( \|\widehat{u} \|^2+ \|\widehat{b} \|^2 \right) \nonumber\\
 	&\leq& \sqrt[p]{2} e_p^{p-1}(k,t)R^2(t). \nonumber
 	\EEQSZ
 	
 	We next put \eqref{one}, \eqref{two} and \eqref{three} together to get,
 	
 	\DEQS
 	\lqq{\frac{d}{dt}e_p^p(k,t)}\\
 	&\leq& -\frac{\min(\nu,\eta)}{4}p|k|^2e_p^p(k,t)+\frac{3|k|}{2}(2\delta)^{3/p}2^{\frac{1}{p}}e_p^{p-1}(k,t)R^2(t)+ \\
 	&& 2\frac{3|k|p}{2}e_p^{p-1}(k,t)h_p(k,t) \\
 	&\leq & -\frac{\min(\nu,\eta)}{4}p|k|^2e_p^p(k,t) + \frac{3|k|}{2}pe_p^{p-1}(k,t)\left(2^\frac{1}{p}(2\delta)^{3/p}R^2(t) +  2h_p(k,t) \right).
 	\EEQS	
 	
 	Once again we consider the set 
 	
 	\DEQS
 	B_{R_1} = \left\lbrace e(k,t):0\leq e(k,t)\leq \frac{R_1(t)}{|k|} \right\rbrace.
 	\EEQS
 	
 	Setting $e(k,t)=e_p(k,t) =\frac{R_1(t)}{|k|} $,   on the boundary such that $ |k|e_p(k,t) = R_1(t)$,
 	
 	\DEQSZ
 	\frac{d}{dt}e_p^p(k,t)&\leq & -\frac{\nu}{4}p|k|^2\frac{R_1^p(t) }{|k|^p} + p\frac{3|k|}{2}\frac{R_1^{p-1}(t) }{|k|^{p-1}} \left( 2^\frac{1}{p}(2\delta)^{3/p}R^2(t) +  2h_p(k,t) \right) \label{attractor2}\\
 	&< & -\frac{\min(\nu,\eta)}{4}p|k|^2\frac{R_1^p(t) }{|k|^p} + p\frac{3|k|}{2}\frac{R_1^{p-1}(t) }{|k|^{p-1}}\frac{\min(\nu,\eta)}{6}R_1(t)=0\nonumber
 	\EEQSZ
 	
 	Here we   used the condition that $\displaystyle 2^\frac{1}{p}(2\delta)^{3/p}R^2(t) +  2h_p(k,t) < \frac{\min(\nu,\eta)}{6}R_1(t) $.
 	Thus, \eqref{attractor2} implies $B_{R_1}$ is an attracting set for $e_p(k,t).$
 	Therefore, if $e_p(k,0)<\frac{R_1(0)}{|k|}$, then $e_p(k,t)<\frac{R_1(t)}{|k|}$ for all $t\in\mathbb{R}^+$.
 \end{proof}

 The following two theorems are the main results of this section, which  
 are direct consequences of Lemma \ref{lemma2.9} and Lemma \ref{lemma2.10}. 

 \begin{Theorem}\label{thm2.7}
 	Let the assumptions of Lemma \ref{lemma2.10} hold. 
 	%
 	%
 	If the weak solution $(u,b)$ of \eqref{mhdmain} satisfies the initial condition
 	
 	\DEQS
 	\sup\limits_{2\leq p<\infty} e_p(k,0)<\frac{R_1(0)}{|k|},
 	\EEQS
 	
 	then for all   $t>0,$
 	
 	\DEQS
 	\sup\limits_{2\leq p<\infty} e_p(k,t)<\frac{R_1(t)}{|k|},
 	\EEQS
 	
 	holds.
 \end{Theorem}

 \begin{Theorem}\label{thm2.8}
 	Suppose the weak solution $(u,b)$ of \eqref{mhdmain} satisfies \eqref{17} 
 	and 
 	$ \sup\limits_{2\leq p<\infty} e_p(k,0)<\dfrac{R_1(0)}{|k|}.$
 	Then for all $T\in\mathbb{R}^+$, we have
 	
 	\DEQSZ
 	\int\limits_0^T\sup\limits_{2\leq p <\infty}e_p(k,t){\rm d} t\leq  \frac{R_2^2(T)}{\min(\nu,\eta)|k|^4}, \label{27*}
 	\EEQSZ
 	
 	and 
 	
 	\DEQSZ
 	R_2(T) := \frac{1}{2}\lk(R_3(T)+\sqrt{4R_1^2(0)+R_3^2(T)}\rk) \label{28}
 	\EEQSZ
 	
 	where
 	
 	\DEQSZ
 	R_3(T) &=& \frac{2R^2(T)}{\min(\nu,\eta)} +\frac{2F_\infty(T)}{\sqrt{\min(\nu,\eta)}} \nonumber, \\
 	F_\infty(T) &=& \sup\limits_{k\in\mathbb{R}\ \{0\}}\left(\int\limits_0^T|\widehat{\chi}_k(\xi)\widehat{f}_1|_{L_\infty}^2+|\widehat{\chi}_k(\xi)\widehat{f}_2|_{L_\infty}^2 dt\right) \label{29}
 	\EEQSZ
 	
 \end{Theorem}

 \begin{proof}[Proof of Theorem \ref{thm2.7}]
 	The proof is very direct. Lemma \ref{lemma2.10} implies that $e_p(k,t)$ is bounded uniformly in $p$. 
 	Then taking the supremum over all $2\leq p<\infty$ concludes the proof.  
 \end{proof}
 
 \begin{proof}[Proof of Theorem \ref{thm2.8}]
 	Recalling the definition of $e_p(k,t)$ from \eqref{ep}, we have 
 	
 	\DEQS
 	e_p^2(k,t)=\left(\int|\widehat{\chi}_k(\xi)\widehat{u}(\xi,t|^p+|\widehat{\chi}_k(\xi)\widehat{b}(\xi,t|^p {\rm d}\xi\right)^\frac{2}{p}.
 	\EEQS
 	
 	Now taking the derivative in time, 
 	
 	\DEQSZ
 	\frac{\partial}{\partial_t}e_p^2(k,t) = \frac{2}{p}\left(\int|\widehat{\chi}_k(\xi)\widehat{u}(\xi,t)|^p+|\widehat{\chi}_k(\xi)\widehat{b}(\xi,t)|^p {\rm d} \xi\right)^{\frac{2}{p}-1}\frac{\partial}{\partial_t}e_p^p(k,t).  \label{now}
 	\EEQSZ
 	
 	We now plug  \eqref{depp} in \eqref{now} to get, 
 	
 	\DEQSZ
 	\lqq{\frac{\partial}{\partial_t}e_p^2(k,t)} \label{dep2}\\
 	&=&\left\lbrack\frac{2}{p}\left(\int|\widehat{\chi}_k(\xi)\widehat{u}(\xi,t)|^p+|\widehat{\chi}_k(\xi)\widehat{b}(\xi,t)|^p {\rm d}\xi\right)^{\frac{2}{p}-1}\right\rbrack \nonumber\\
 	&&\times\left\lbrack-\nu\int p|\xi|^2|\widehat{\chi}_k(\xi) \widehat{u}(\xi,t)|^p{\rm d}\xi-\eta\int p|\xi|^2|\widehat{\chi}_k(\xi) \widehat{b}(\xi,t)|^p{\rm d}\xi\right. \nonumber\\
 	&&+\Re\int ip|\widehat{\chi}_k(\xi) \widehat{u}(\xi,t)|^{p-2}\overline{\widehat{\chi}_k(\xi) \widehat{u}(\xi,t) }\widehat{\chi}_k(\xi)\Pi_\xi\int\widehat{u}(\xi-\zeta) \zeta\widehat{u}(\zeta){\rm d}\zeta {\rm d}\xi \nonumber\\
 	&&+\Re\int ip|\widehat{\chi}_k(\xi) \widehat{u}(\xi,t)|^{p-2}\overline{\widehat{\chi}_k(\xi) \widehat{u}(\xi,t) }\widehat{\chi}_k(\xi)\Pi_\xi\int\widehat{b}(\xi-\zeta) \zeta\widehat{b}(\zeta){\rm d}\zeta {\rm d}\xi \nonumber\\
 	&&+\Re\int ip|\widehat{\chi}_k(\xi) \widehat{b}(\xi,t)|^{p-2}\overline{\widehat{\chi}_k(\xi) \widehat{b}(\xi,t) }\widehat{\chi}_k(\xi)\Pi_\xi\int\widehat{b}(\xi-\zeta) \zeta\widehat{u}(\zeta){\rm d}\zeta {\rm d}\xi \nonumber\\
 	&&+\Re\int ip|\widehat{\chi}_k(\xi) \widehat{b}(\xi,t)|^{p-2}\overline{\widehat{\chi}_k(\xi) \widehat{b}(\xi,t) }\widehat{\chi}_k(\xi)\Pi_\xi\int\widehat{u}(\xi-\zeta) \zeta\widehat{b}(\zeta){\rm d}\zeta {\rm d}\xi \nonumber\\
 	&&+\Re\int p|\widehat{\chi}_k(\xi) \widehat{u}(\xi,t)|^{p-2}\overline{\widehat{\chi}_k(\xi) \widehat{u}(\xi,t) }\widehat{\chi}_k(\xi)\widehat{f}_1(\xi,t){\rm d}\xi \nonumber\\
 	&&\left.+\Re\int p|\widehat{\chi}_k(\xi) \widehat{b}(\xi,t)|^{p-2}\overline{\widehat{\chi}_k(\xi) \widehat{b}(\xi,t) }\widehat{\chi}_k(\xi)\widehat{f}_2(\xi,t){\rm d}\xi\right\rbrack \nonumber  .
 	\EEQSZ
 	
 	For the sake of calculation simplicity, we split the RHS of \eqref{dep2} in to the following integrals. 
 	
 	\DEQS
 	I_1 &:= &-2\nu\left(\int|\widehat{\chi}_k(\xi)\widehat{u}(\xi,t)|^p+|\widehat{\chi}_k(\xi)\widehat{b}(\xi,t)|^p{\rm d}\xi \right)^{\frac{2}{p}-1}\\
 	&&\times\left(\int |\xi|^2|\widehat{\chi}_k(\xi) \widehat{u}(\xi,t)|^p +|\xi|^2|\widehat{\chi}_k(\xi) \widehat{b}(\xi,t)|^p{\rm d}\xi  \right),\\
 	I_2 &:= & 2\left(\int|\widehat{\chi}_k(\xi)\widehat{u}(\xi,t)|^p+|\widehat{\chi}_k(\xi)\widehat{b}(\xi,t)|^p{\rm d}\xi \right)^{\frac{2}{p}-1}\\
 	&&\times\Re\int \left( i|\widehat{\chi}_k(\xi) \widehat{u}(\xi,t)|^{p-2}\overline{\widehat{\chi}_k(\xi) \widehat{u}(\xi,t) }\widehat{\chi}_k(\xi)\Pi_\xi\int\widehat{u}(\xi-\zeta) \zeta\widehat{u}(\zeta){\rm d}\zeta \right){\rm d}\xi,\\
 	I_3 &:= & 2\left(\int|\widehat{\chi}_k(\xi)\widehat{u}(\xi,t)|^p+|\widehat{\chi}_k(\xi)\widehat{b}(\xi,t)|^p{\rm d}\xi \right)^{\frac{2}{p}-1}\\
 	&&\times\Re\int\left(  i|\widehat{\chi}_k(\xi) \widehat{u}(\xi,t)|^{p-2}\overline{\widehat{\chi}_k(\xi) \widehat{u}(\xi,t) }\widehat{\chi}_k(\xi)\Pi_\xi\int\widehat{b}(\xi-\zeta) \zeta\widehat{b}(\zeta){\rm d}\zeta\right) {\rm d}\xi,\\
 	I_4 &:= & 2\left(\int|\widehat{\chi}_k(\xi)\widehat{u}(\xi,t)|^p+|\widehat{\chi}_k(\xi)\widehat{b}(\xi,t)|^p{\rm d}\xi \right)^{\frac{2}{p}-1}\\
 	&&\times\Re\int\left(  i|\widehat{\chi}_k(\xi) \widehat{b}(\xi,t)|^{p-2}\overline{\widehat{\chi}_k(\xi) \widehat{b}(\xi,t) }\widehat{\chi}_k(\xi)\Pi_\xi\int\widehat{b}(\xi-\zeta) \zeta\widehat{u}(\zeta){\rm d}\zeta\right) {\rm d}\xi,\\
 	I_5 &:= & 2\left(\int|\widehat{\chi}_k(\xi)\widehat{u}(\xi,t)|^p+|\widehat{\chi}_k(\xi)\widehat{b}(\xi,t)|^p{\rm d}\xi \right)^{\frac{2}{p}-1}\\
 	&&\times\Re\int\left(  i|\widehat{\chi}_k(\xi) \widehat{b}(\xi,t)|^{p-2}\overline{\widehat{\chi}_k(\xi) \widehat{b}(\xi,t) }\widehat{\chi}_k(\xi)\Pi_\xi\int\widehat{u}(\xi-\zeta) \zeta\widehat{b}(\zeta){\rm d}\zeta\right) {\rm d}\xi, \\ 
 	I_6 &:= & 2\left(\int|\widehat{\chi}_k(\xi)\widehat{u}(\xi,t)|^p+|\widehat{\chi}_k(\xi)\widehat{b}(\xi,t)|^p{\rm d}\xi \right)^{\frac{2}{p}-1} \\
 	&&\times\Re\int\left( |\widehat{\chi}_k(\xi) \widehat{u}(\xi,t)|^{p-2}\overline{\widehat{\chi}_k(\xi) \widehat{u}(\xi,t) }\widehat{\chi}_k(\xi)\widehat{f}_1(\xi,t) \right){\rm d}\xi , \\
 	I_7 &:= & 2\left(\int|\widehat{\chi}_k(\xi)\widehat{u}(\xi,t)|^p+|\widehat{\chi}_k(\xi)\widehat{b}(\xi,t)|^p{\rm d}\xi \right)^{\frac{2}{p}-1} \\
 	&&\times\Re\int\left( |\widehat{\chi}_k(\xi) \widehat{b}(\xi,t)|^{p-2}\overline{\widehat{\chi}_k(\xi) \widehat{b}(\xi,t) }\widehat{\chi}_k(\xi)\widehat{f}_2(\xi,t)\right){\rm d}\xi.
 	\EEQS
 	
 	We now proceed to  estimating each of these integrals $(I_1) - (I_7)$.
 	
 	\DEQSZ
 	I_1 &=& -2\nu\left(\int|\widehat{\chi}_k(\xi)\widehat{u}(\xi,t)|^p+|\widehat{\chi}_k(\xi)\widehat{b}(\xi,t)|^p{\rm d}\xi \right)^{\frac{2}{p}-1} \label{I12est3}r\\
 	&&\left(\int |\xi|^2|\widehat{\chi}_k(\xi) \widehat{u}(\xi,t)|^p +|\xi|^2|\widehat{\chi}_k(\xi) \widehat{b}(\xi,t)|^p{\rm d}\xi  \right)^\frac{2}{p} \nonumber\\
 	& =& -2\nu\left( \frac{ \int|\widehat{\chi}_k(\xi)\widehat{u}(\xi,t)|^p+|\widehat{\chi}_k(\xi)\widehat{b}(\xi,t)|^p{\rm d}\xi }{\int |\xi|^2(|\widehat{\chi}_k(\xi) \widehat{u}(\xi,t)|^p +|\widehat{\chi}_k(\xi) \widehat{b}(\xi,t)|^p){\rm d}\xi} \right)^{\frac{2}{p}-1} \nonumber\\
 	&&\left(\int |\xi|^2(|\widehat{\chi}_k(\xi) \widehat{u}(\xi,t)|^p +|\widehat{\chi}_k(\xi) \widehat{b}(\xi,t)|^p){\rm d}\xi  \right)^{2/p} \nonumber.
 	\EEQSZ
 	
 	\DEQSZ
 	|I_2| &=& \left|  2\left(\int|\widehat{\chi}_k(\xi)\widehat{u}(\xi,t)|^p+|\widehat{\chi}_k(\xi)\widehat{b}(\xi,t)|^p{\rm d}\xi \right)^{\frac{2}{p}-1} \right.\label{I2est3}\\
 	&&\times\left.  \int\left( i|\widehat{\chi}_k(\xi) \widehat{u}(\xi,t)|^{p-2}\overline{\widehat{\chi}_k(\xi) \widehat{u}(\xi,t) }\widehat{\chi}_k(\xi)\Pi_\xi\int\widehat{u}(\xi-\zeta) \zeta\widehat{u}(\zeta){\rm d}\zeta \right){\rm d}\xi\right|\nonumber\\
 	&\leq& 2 \left(\int|\widehat{\chi}_k(\xi)\widehat{u}(\xi,t)|^p+|\widehat{\chi}_k(\xi)\widehat{b}(\xi,t)|^p{\rm d}\xi \right)^{\frac{2}{p}-1} \left(\int|\widehat{\chi}_k(\xi)\widehat{u}(\xi,t)|^p {\rm d}\xi\right)^{\frac{p-2}{p} }\nonumber\\
 	&& \times \left(\int|\widehat{\chi}_k(\xi)\widehat{u}(\xi,t)|^p {\rm d}\xi\right)^{\frac{1}{p}} \left(\int|\widehat{\chi}_k(\xi)|^p{\rm d}\xi \right)^{\frac{1}{p}}   \| \Pi_\xi\int\widehat{u}(\xi-\zeta) \zeta\widehat{u}(\zeta){\rm d}\zeta  \|_{L^\infty}. \nonumber
 	\EEQSZ
 	
 	Here we repeatedly used Holder's inequality.
 	Similar calculations give us
 	
 	\DEQSZ
 	|I_3| &\leq& 2 \left(\int|\widehat{\chi}_k(\xi)\widehat{u}(\xi,t)|^p+|\widehat{\chi}_k(\xi)\widehat{b}(\xi,t)|^p{\rm d}\xi \right)^{\frac{2}{p}-1}\left(\int|\widehat{\chi}_k(\xi)\widehat{u}(\xi,t)|^p {\rm d}\xi\right)^{\frac{p-2}{p} }\label{I3est3}\\
 	&&\times\left(\int|\widehat{\chi}_k(\xi)\widehat{u}(\xi,t)|^p {\rm d}\xi\right)^{\frac{1}{p}} \left(\int|\widehat{\chi}_k(\xi)|^p{\rm d}\xi \right)^{\frac{1}{p}}   \| \Pi_\xi\int\widehat{b}(\xi-\zeta) \zeta\widehat{b}(\zeta){\rm d}\zeta  \|_{L^\infty} \nonumber,
 	\EEQSZ
 	
 	\DEQSZ
 	|I_4| &\leq& 2 \left(\int|\widehat{\chi}_k(\xi)\widehat{u}(\xi,t)|^p+|\widehat{\chi}_k(\xi)\widehat{b}(\xi,t)|^p{\rm d}\xi \right)^{\frac{2}{p}-1}\left(\int|\widehat{\chi}_k(\xi)\widehat{b}(\xi,t)|^p {\rm d}\xi\right)^{\frac{p-2}{p} }\label{I4est3}\\
 	&&\times\left(\int|\widehat{\chi}_k(\xi)\widehat{b}(\xi,t)|^p {\rm d}\xi\right)^{\frac{1}{p}} \left(\int|\widehat{\chi}_k(\xi)|^p{\rm d}\xi \right)^{\frac{1}{p}}   \| \Pi_\xi\int\widehat{b}(\xi-\zeta) \zeta\widehat{u}(\zeta){\rm d}\zeta  \|_{L^\infty}\nonumber,
 	\EEQSZ
 	
 	\DEQSZ
 	|I_5|  &\leq& 2 \left(\int|\widehat{\chi}_k(\xi)\widehat{u}(\xi,t)|^p+|\widehat{\chi}_k(\xi)\widehat{b}(\xi,t)|^p{\rm d}\xi \right)^{\frac{2}{p}-1}\left(\int|\widehat{\chi}_k(\xi)\widehat{b}(\xi,t)|^p {\rm d}\xi\right)^{\frac{p-2}{p} }\label{I5est3}\\
 	&&\times\left(\int|\widehat{\chi}_k(\xi)\widehat{b}(\xi,t)|^p {\rm d}\xi\right)^{\frac{1}{p}} \left(\int|\widehat{\chi}_k(\xi)|^p{\rm d}\xi \right)^{\frac{1}{p}}   \| \Pi_\xi\int\widehat{u}(\xi-\zeta) \zeta\widehat{b}(\zeta){\rm d}\zeta  \|_{L^\infty} \nonumber.
 	\EEQSZ
 	
 	For integrals involving  the inhomogeneous forces,
 	
 	\DEQSZ
 	|I_6|&\leq&2\left(\int|\widehat{\chi}_k(\xi)\widehat{u}(\xi,t)|^p+|\widehat{\chi}_k(\xi)\widehat{b}(\xi,t)|^p{\rm d}\xi \right)^{\frac{2}{p}-1} \label{I6est3} \\
 	&&\times\left|\int |\widehat{\chi}_k(\xi) \widehat{u}(\xi,t)|^{p-2}\overline{\widehat{\chi}_k(\xi) \widehat{u}(\xi,t) }\widehat{\chi}_k(\xi)\widehat{f}_1(\xi,t){\rm d}\xi   \right|\nonumber\\
 	&\leq& 2 \left(\int|\widehat{\chi}_k(\xi)\widehat{u}(\xi,t)|^p+|\widehat{\chi}_k(\xi)\widehat{b}(\xi,t)|^p{\rm d}\xi \right)^{\frac{2}{p}-1} \left(\int|\widehat{\chi}_k(\xi)\widehat{u}(\xi,t)|^p {\rm d}\xi\right)^{\frac{p-2}{p} }\nonumber\\
 	&& \times\left(\int|\widehat{\chi}_k(\xi)\widehat{u}(\xi,t)|^p {\rm d}\xi\right)^{\frac{1}{p} }\left(\int|\widehat{\chi}_k(\xi)\widehat{f}_1(\xi,t)|^p{\rm d}\xi \right)^{\frac{1}{p}} \nonumber .
 	\EEQSZ
 	
 	Similarly,
 	
 	\DEQSZ
 	|I_7| &\leq&2 \left(\int|\widehat{\chi}_k(\xi)\widehat{u}(\xi,t)|^p+|\widehat{\chi}_k(\xi)\widehat{b}(\xi,t)|^p{\rm d}\xi \right)^{\frac{2}{p}-1} \left(\int|\widehat{\chi}_k(\xi)\widehat{b}(\xi,t)|^p {\rm d}\xi\right)^{\frac{p-2}{p} }\label{I7est3}\\
 	&&\times \left(\int|\widehat{\chi}_k(\xi)\widehat{b}(\xi,t)|^p {\rm d}\xi\right)^{\frac{1}{p} }\left(\int|\widehat{\chi}_k(\xi)\widehat{f}_2(\xi,t)|^p{\rm d}\xi \right)^{\frac{1}{p}} .\nonumber
 	\EEQSZ  
 	
 	Now taking the time integral of \eqref{dep2} over the interval $[0,T]$ we get
 	
 	\DEQS
 	e_p^2(k,T) - e_p^2(k,0) = \int\limits_0^T \sum\limits_{j=1}^7I_j{\rm d} t.
 	\EEQS
 	
 	Then it follows from \eqref{I12est3} that,
 	
 	\DEQSZ
 	\lqq{2\min(\nu,\eta) \bigints\limits_0^T\Big\{\left(\frac{\int |\xi|^2(|\widehat{\chi}_k(\xi) \widehat{u}(\xi,t)|^p +|\widehat{\chi}_k(\xi) \widehat{b}(\xi,t)|^p){\rm d}\xi }{ \int|\widehat{\chi}_k(\xi)\widehat{u}(\xi,t)|^p+|\widehat{\chi}_k(\xi)\widehat{b}(\xi,t)|^p{\rm d}\xi  } \right)^{1-\frac{2}{p}}}\label{finI1} \\
 	&&\times\left(\int |\xi|^2(|\widehat{\chi}_k(\xi) \widehat{u}(\xi,t)|^p +|\widehat{\chi}_k(\xi) \widehat{b}(\xi,t)|^p){\rm d}\xi  \right)^{2/p}\Big\}{\rm d} t\nonumber\\
 	&\leq& e_p^2(k,0) - e_p^2(k,T) + \sum\limits_{j=2}^7\int\limits_0^T |I_j|{\rm d} t. \nonumber
 	\EEQSZ
 	
 	Once again making use of the Young's inequality gives,
 	
 	\DEQS
 	\| \Pi_\xi\int\widehat{u}(\xi-\zeta) \zeta\widehat{b}(\zeta)d\zeta  \|_{L^\infty} \leq  \|\widehat{u}(\cdot,t) \|_{L^2}  \|\xi\widehat{b}(\cdot,t) \|_{L^2}.
 	\EEQS
 	
 	Therefore,
 	
 	\DEQSZ
 	\lqq{\int\limits_0^T|I_2|{\rm d} t }\label{testI2}\\
 	&\leq&  2\int\limits_0^T\left( \left(\int|\widehat{\chi}_k(\xi)\widehat{u}(\xi,t)|^p+|\widehat{\chi}_k(\xi)\widehat{b}(\xi,t)|^p{\rm d}\xi \right)^{\frac{2}{p}-1} \left(\int|\widehat{\chi}_k(\xi)\widehat{u}(\xi,t)|^p {\rm d}\xi\right)^{\frac{p-2}{p} }\right. \nonumber \\
 	&&\left(\int|\widehat{\chi}_k(\xi)\widehat{u}(\xi,t)|^p {\rm d}\xi\right)^{\frac{1}{p}}\left.\quad \left(\int|\widehat{\chi}_k(\xi)|^p{\rm d}\xi \right)^{\frac{1}{p}}   \| \Pi_\xi\int\widehat{u}(\xi-\zeta) \zeta\widehat{u}(\zeta){\rm d}\zeta  \|_{L^\infty}\right){\rm d} t \nonumber\\
 	&&\leq   2\int\limits_0^T\left( \left(\int|\widehat{\chi}_k(\xi)\widehat{u}(\xi,t)|^p+|\widehat{\chi}_k(\xi)\widehat{b}(\xi,t)|^p{\rm d}\xi \right)^{\frac{2}{p}-1} \left(\int|\widehat{\chi}_k(\xi)\widehat{u}(\xi,t)|^p {\rm d}\xi\right)^{\frac{p-1}{p} }\right.\nonumber\\
 	&&\left.\quad \left(\int|\widehat{\chi}_k(\xi)|^p{\rm d}\xi \right)^{\frac{1}{p}}  \|\widehat{u}(\cdot,t) \|_{L^2}  \|\xi\widehat{u}(\cdot,t) \|_{L^2} \right){\rm d} t \nonumber.
 	\EEQSZ
 	
 	Thus, similar computations yield,
 	
 	\DEQSZ
 	\int\limits_0^T |I_3|{\rm d} t&\leq& 2\int\limits_0^T\left( \left(\int|\widehat{\chi}_k(\xi)\widehat{u}(\xi,t)|^p+|\widehat{\chi}_k(\xi)\widehat{b}(\xi,t)|^p{\rm d}\xi \right)^{\frac{2}{p}-1} \right.  \label{testI3}\\
 	&&\left.\left(\int|\widehat{\chi}_k(\xi)\widehat{u}(\xi,t)|^p {\rm d}\xi\right)^{\frac{p-1}{p} } \left(\int|\widehat{\chi}_k(\xi)|^p{\rm d}\xi \right)^{\frac{1}{p}}  \|\widehat{b}(\cdot,t) \|_{L^2}  \|\xi\widehat{b}(\cdot,t) \|_{L^2} \right){\rm d} t \nonumber,\\
 	\int\limits_0^T|I_4| {\rm d}t &\leq &  2\int\limits_0^T\left( \left(\int|\widehat{\chi}_k(\xi)\widehat{u}(\xi,t)|^p+|\widehat{\chi}_k(\xi)\widehat{b}(\xi,t)|^p{\rm d}\xi \right)^{\frac{2}{p}-1} \right.\label{testI4} \\
 	&&\left.\left(\int|\widehat{\chi}_k(\xi)\widehat{b}(\xi,t)|^p{\rm d}\xi \right)^{\frac{p-1}{p} } \left(\int|\widehat{\chi}_k(\xi)|^p{\rm d}\xi \right)^{\frac{1}{p}}  \|\widehat{u}(\cdot,t) \|_{L^2}  \|\xi\widehat{b}(\cdot,t) \|_{L^2} \right){\rm d} t  \nonumber,\\
 	\int\limits_0^T|I_5|{\rm d} t &\leq &  2\int\limits_0^T\left( \left(\int|\widehat{\chi}_k(\xi)\widehat{u}(\xi,t)|^p+|\widehat{\chi}_k(\xi)\widehat{b}(\xi,t)|^p{\rm d}\xi \right)^{\frac{2}{p}-1} \right.\label{testI5} \\
 	&&\left.\left(\int|\widehat{\chi}_k(\xi)\widehat{b}(\xi,t)|^p {\rm d}\xi\right)^{\frac{p-1}{p} } \left(\int|\widehat{\chi}_k(\xi)|^p{\rm d}\xi \right)^{\frac{1}{p}}  \|\widehat{b}(\cdot,t) \|_{L^2}  \|\xi\widehat{u}(\cdot,t) \|_{L^2} \right){\rm d} t \nonumber,\\
 	\int\limits_0^T|I_6|{\rm d} t&\leq &  2\int\limits_0^T{ \left(\int|\widehat{\chi}_k(\xi)\widehat{u}(\xi,t)|^p+|\widehat{\chi}_k(\xi)\widehat{b}(\xi,t)|^p{\rm d}\xi \right)^{\frac{2}{p}-1}}  \label{testI6}\\
 	&& \left(\int|\widehat{\chi}_k(\xi)\widehat{u}(\xi,t)|^p {\rm d}\xi\right)^{\frac{p-1}{p} } \left(\int|\widehat{\chi}_k(\xi)\widehat{f}_1(\xi,t)|^p{\rm d}\xi \right)^{\frac{1}{p}}{\rm d} t \nonumber,\\
 	\int\limits_0^T|I_7|{\rm d} t&\leq &  2\int\limits_0^T{ \left(\int|\widehat{\chi}_k(\xi)\widehat{u}(\xi,t)|^p+|\widehat{\chi}_k(\xi)\widehat{b}(\xi,t)|^p{\rm d}\xi \right)^{\frac{2}{p}-1} }    \label{testI7}\\
 	&& \left(\int|\widehat{\chi}_k(\xi)\widehat{b}(\xi,t)|^p {\rm d}\xi\right)^{\frac{p-1}{p} } \left(\int|\widehat{\chi}_k(\xi)\widehat{f}_2(\xi,t)|^p{\rm d}\xi \right)^{\frac{1}{p}} {\rm d} t \nonumber.
 	\EEQSZ
 	
 	Now putting estimates \eqref{testI2}-\eqref{testI5} together  we get,
 	
 	\DEQSZ
 	\lqq{\int\limits_0^T|I_2|{\rm d} t +\int\limits_0^T|I_3|{\rm d} t+\int\limits_0^T|I_4|{\rm d} t+\int\limits_0^T|I_5|{\rm d} t} \label{finI25}\\
 	&\leq& 2 \left(\int|\widehat{\chi}_k(\xi)|^p{\rm d}\xi \right)^{\frac{1}{p}}\bigintss_0^T\Big\{ \left(\int|\widehat{\chi}_k(\xi)\widehat{u}(\xi,t)|^p+|\widehat{\chi}_k(\xi)\widehat{b}(\xi,t)|^p{\rm d}\xi \right)^{\frac{2}{p}-1}  \nonumber\\
 	&&\times \left[ \left(\int|\widehat{\chi}_k(\xi)\widehat{u}(\xi,t)|^p {\rm d}\xi\right)^{\frac{p-1}{p} } \left( \|\widehat{u}(\cdot,t) \|_{L^2}  \|\xi\widehat{u}(\cdot,t) \|_{L^2} +  \|\widehat{b}(\cdot,t) \|_{L^2}  \|\xi\widehat{b}(\cdot,t) \|_{L^2} \right) \right.\nonumber\\
 	&&+\left.\left(\int|\widehat{\chi}_k(\xi)\widehat{b}(\xi,t)|^p \right)^{\frac{p-1}{p} }\left(  \|\widehat{u}(\cdot,t) \|_{L^2}  \|\xi\widehat{b}(\cdot,t) \|_{L^2} +  \|\widehat{b}(\cdot,t) \|_{L^2} \|\xi\widehat{u}(\cdot,t) \|_{L^2} \right) \right] \Big\}{\rm d} t \nonumber\\
 	&\leq& 2 \left(\int|\widehat{\chi}_k(\xi)|^p{\rm d}\xi \right)^{\frac{1}{p}}\bigintss_0^T\Big\{\left(\int|\widehat{\chi}_k(\xi)\widehat{u}(\xi,t)|^p+|\widehat{\chi}_k(\xi)\widehat{b}(\xi,t)|^p{\rm d}\xi \right)^{\frac{2}{p}-1}  \nonumber\\
 	&&\times \left[ \left(\int|\widehat{\chi}_k(\xi)\widehat{u}(\xi,t)|^p + |\widehat{\chi}_k(\xi)\widehat{b}(\xi,t)|^p {\rm d}\xi\right)^{\frac{p-1}{p} } \right. \nonumber\\
 	&&\left.\times  \left( \|\xi\widehat{u}(\cdot,t) \|_{L^2} + \|\xi\widehat{b}(\cdot,t) \|_{L^2}  \right) \left( \|\widehat{u}(\cdot,t) \|_{L^2}+ \|\widehat{b}(\cdot,t) \|_{L^2} \right)\right]\Big\}{\rm d} t \nonumber\\
 	&\leq&  2(2\delta)^{\frac{3}{p}} \int\limits_0^T{\left( \int |\widehat{\chi}_k\widehat{u}|^p+|\widehat{\chi}_k\widehat{b}|^p{\rm d} \xi\right)^\frac{1}{p}\left( \|\widehat{u} \|+ \|\widehat{b} \| \right)\left(  \|\xi\widehat{u} \|+ \|\xi\widehat{b} \|\right) }{\rm d} t \nonumber \\
 	&\leq&   2(2\delta)^{\frac{3}{p}} \left( \int\limits_0^T\left( \int |\widehat{\chi}_k\widehat{u}|^p+|\widehat{\chi}_k\widehat{b}|^p{\rm d} \xi\right)^\frac{2}{p}\right)^\frac{1}{2}\sup_{0\leq t\leq T}\left( \|\widehat{u} \|+ \|\widehat{b} \| \right)\nonumber\\
 	&&\left(\int\limits_0^T\left(  \|\xi\widehat{u} \|+ \|\xi\widehat{b} \|\right)^2{\rm d} t\right)^\frac{1}{2} \nonumber\\
 	&\leq &  2(2\delta)^{\frac{3}{p}} {\rm supp }_{0\leq t\leq T}\left( \|\widehat{u} \|+ \|\widehat{b} \| \right)\left( \int\limits_0^T\left( \int |\widehat{\chi}_k\widehat{u}|^p+|\widehat{\chi}_k\widehat{b}|^p{\rm d} \xi\right)^\frac{2}{p}\right)^\frac{1}{2} \nonumber\\
 	&&\left(\int\limits_0^T\left(  \|\nabla u \|+ \|\nabla b \|\right)^2{\rm d} t\right)^\frac{1}{2} \nonumber \\
 	&\leq & 2(2\delta)^\frac{3}{p}R^2(T)\left( \int\limits_0^T\left( \int |\widehat{\chi}_k\widehat{u}|^p+|\widehat{\chi}_k\widehat{b}|^p{\rm d} \xi\right)^\frac{2}{p}\right)^\frac{1}{2} \nonumber.
 	\EEQSZ
 	
 	And, from \eqref{testI6} and \eqref{testI7} we have,
 	
 	\DEQSZ
 	\lqq{\int\limits_0^T|I_6|{\rm d} t +\int\limits_0^T|I_7|{\rm d} t} \nonumber\\
 	&\leq & 2\int\limits_0^T\left\lbrace\left(\int|\widehat{\chi}_k(\xi)\widehat{u}(\xi,t)|^p+|\widehat{\chi}_k(\xi)\widehat{b}(\xi,t)|^p{\rm d}\xi \right)^{\frac{2}{p}-1}\right.\nonumber\\ 
 	&& \left[\left(\int|\widehat{\chi}_k(\xi)\widehat{u}(\xi,t)|^p {\rm d}\xi\right)^{\frac{p-1}{p} } \left(\int|\widehat{\chi}_k(\xi)\widehat{f}_1(\xi,t)|^p{\rm d}\xi \right)^{\frac{1}{p}} \right.\nonumber \\
 	&& \left.+\left.  \left(\int|\widehat{\chi}_k(\xi)\widehat{b}(\xi,t)|^p{\rm d}\xi \right)^{\frac{p-1}{p} } \left(\int|\widehat{\chi}_k(\xi)\widehat{f}_2(\xi,t)|^p{\rm d}\xi \right)^{\frac{1}{p}} \right]\right\rbrace{\rm d} t \nonumber\\
 	&\leq & 2\int\limits_0^T\left\lbrace \left(\int|\widehat{\chi}_k(\xi)\widehat{u}(\xi,t)|^p+|\widehat{\chi}_k(\xi)\widehat{b}(\xi,t)|^p{\rm d}\xi \right)^{\frac{2}{p}-1} \right. \nonumber\\
 	&& \left(\left(\int|\widehat{\chi}_k(\xi)\widehat{u}(\xi,t)|^p {\rm d}\xi\right)^{\frac{p-1}{p} } +\left(\int|\widehat{\chi}_k(\xi)\widehat{b}(\xi,t)|^p {\rm d}\xi\right)^{\frac{p-1}{p} }\right) \nonumber\\
 	&& \left. \left(\left(\int|\widehat{\chi}_k(\xi)\widehat{f}_2(\xi,t)|^p{\rm d}\xi \right)^{\frac{1}{p}}+\left(\int|\widehat{\chi}_k(\xi)\widehat{f}_2(\xi,t)|^p{\rm d}\xi \right)^{\frac{1}{p}} \right) \right\rbrace{\rm d} t\nonumber\\
 	&\leq &  2\int\limits_0^T\left\lbrace \left(\int|\widehat{\chi}_k(\xi)\widehat{u}(\xi,t)|^p+|\widehat{\chi}_k(\xi)\widehat{b}(\xi,t)|^p{\rm d}\xi \right)^{\frac{2}{p}-1}  \right. \nonumber\\
 	&& \left( \int|\widehat{\chi}_k(\xi)\widehat{u}(\xi,t)|^p+|\widehat{\chi}_k(\xi)\widehat{b}(\xi,t)|^p{\rm d}\xi\right)^{\frac{p-1}{p} } \nonumber \\
 	&& \left. \left( \int|\widehat{\chi}_k(\xi)\widehat{f}_2(\xi,t)|^p+|\widehat{\chi}_k(\xi)\widehat{f}_2(\xi,t)|^p{\rm d}\xi \right)^{\frac{1}{p}}\right\rbrace{\rm d} t \nonumber\\
 	&\leq&  2\left(\int\limits_0^T\left(\int|\widehat{\chi}_k(\xi)\widehat{u}(\xi,t)|^p+|\widehat{\chi}_k(\xi)\widehat{b}(\xi,t)|^p{\rm d}\xi \right)^{\frac{2}{p}}{\rm d} t\right)^{\frac{1}{2}}\nonumber\\
 	&&\left(\int\limits_0^T\left(\int|\widehat{\chi}_k(\xi)\widehat{f}_1(\xi,t)|^p+|\widehat{\chi}_k(\xi)\widehat{f}_2(\xi,t)|^p{\rm d}\xi \right)^{\frac{2}{p}}{\rm d} t\right)^{\frac{1}{2}}. \label{finI67}
 	\EEQSZ
 	
 	Therefore putting \eqref{finI1}, \eqref{finI25} and \eqref{finI67} together and using the fact that $|\xi|\geq \dfrac{|k|}{2}$ in the support of $\widehat{\chi}_k$ gives,
 	
 	\DEQSZ
 	\lqq{\frac{1}{2}\min(\nu,\eta)|k|^2 \int\limits_0^T\left(\int|\widehat{\chi}_k(\xi)\widehat{u}(\xi,t)|^p+|\widehat{\chi}_k(\xi)\widehat{b}(\xi,t)|^p{\rm d}\xi \right)^{\frac{2}{p}}{\rm d} t}\label{ttt}\\
 	&\leq&2\left(\int\limits_0^T\left(\int|\widehat{\chi}_k(\xi)\widehat{u}(\xi,t)|^p+|\widehat{\chi}_k(\xi)\widehat{b}(\xi,t)|^p{\rm d}\xi \right)^{\frac{2}{p}}{\rm d} t\right)^{\frac{1}{2}}\left\lbrack 2(2\delta)^{\frac{3}{p}}R^2(t) +\qquad \right.\nonumber\\
 	&&\left.\left(\int\limits_0^T\left(\int|\widehat{\chi}_k(\xi)\widehat{f}_1(\xi,t)|^p+|\widehat{\chi}_k(\xi)\widehat{f}_2(\xi,t)|^p{\rm d}\xi \right)^{\frac{2}{p}}{\rm d} t\right)^{\frac{1}{2}}  \right] + e_p^2(k,0)- e_p^2(k,T) \nonumber.
 	\EEQSZ
 	
 	Now multiplying \eqref{ttt} by $|k|^2$,
 	
 	\DEQSZ
 	\lqq{\min(\nu,\eta)|k|^4 \int\limits_0^T\left(\int|\widehat{\chi}_k(\xi)\widehat{u}(\xi,t)|^p+|\widehat{\chi}_k(\xi)\widehat{b}(\xi,t)|^p{\rm d}\xi \right)^{\frac{2}{p}}{\rm d} t} \label{tofinal} \\
 	&\leq & {4}\left(\int\limits_0^T|k|^4\left(\int|\widehat{\chi}_k(\xi)\widehat{u}(\xi,t)|^p+|\widehat{\chi}_k(\xi)\widehat{b}(\xi,t)|^p{\rm d}\xi \right)^{\frac{2}{p}}{\rm d} t\right)^{\frac{1}{2}} \nonumber\\
 	&&\left[ (2\delta)^{\frac{3}{p}}R^2(t) +\left(\int\limits_0^T\left(\int|\widehat{\chi}_k(\xi)\widehat{f}_1(\xi,t)|^p+|\widehat{\chi}_k(\xi)\widehat{f}_2(\xi,t)|^p{\rm d}\xi \right)^{\frac{2}{p}}{\rm d} t\right)^{\frac{1}{2}}  \right] \\
 	&&+ {} \, |k|^2\left[e_p^2(k,0)- e_p^2(k,T)\right]. \nonumber
 	\EEQSZ
 	
 	Define,
 	
 	\DEQSZ
 	I_p^2(k,T)&=&\bigintsss\limits_0^T|k|^4\left(\int|\widehat{\chi}_k(\xi)\widehat{u}(\xi,t)|^p+|\widehat{\chi}_k(\xi)\widehat{b}(\xi,t)|^p{\rm d}\xi \right)^{\frac{2}{p}}{\rm d} t \label{ip2},\\
 	\mathbb{F}_p(T) &=& \left(\int\limits_0^T\left(\int|\widehat{\chi}_k(\xi)\widehat{f}_1(\xi,t)|^p+|\widehat{\chi}_k(\xi)\widehat{f}_2(\xi,t)|^p{\rm d}\xi \right)^{\frac{2}{p}}{\rm d} t\right)^{\frac{1}{2}} . \label{fp2}
 	\EEQSZ
 	
 	We now put \eqref{tofinal}-\eqref{fp2} together, use the assumption $e_p(k,0)\leq\frac{R_1(0)}{|k|}$ and rearrange terms to get
 	
 	\DEQSZ
 	\min(\nu,\eta) I_p^2(k,T)- {4}\left[(2\delta)^{\frac{3}{p}}R^2(T)+\mathbb{F}_p(T) \right]I_p(k,T)-  R_1^2(0)\leq 0 \label{crucial}.
 	\EEQSZ
 	
 	Observe that \eqref{crucial} is quadratic in $I_p$. Solving the associated quadratic equation yields
 	
 	\DEQS
 	\frac{4\left[(2\delta)^{\frac{3}{p}}R^2(T)+\mathbb{F}_p(T) \right] \pm \sqrt{\left(4\left[(2\delta)^{\frac{3}{p}}R^2(T)+\mathbb{F}_p(T) \right]\right)^2 + 4\min(\nu,\eta) R_1^2(0)  } }{2\min(\nu,\eta)} .
 	\EEQS
 	
 	Elementary mathematics tells us that $ I_p(k,t)$ cannot exceed the largest positive root of the associated quadratic equation, which is
 	
 	\DEQS
 	\frac{4\left[(2\delta)^{\frac{3}{p}}R^2(T)+\mathbb{F}_p(T) \right] + \sqrt{\left(4\left[(2\delta)^{\frac{3}{p}}R^2(T)+\mathbb{F}_p(T) \right]\right)^2 + 4\min(\nu,\eta) R_1^2(0)  } }{2\min(\nu,\eta)} .
 	\EEQS
 	
 	Now set,
 	
 	\DEQS
 	R_{3,p}(T):=4\lk[(2\delta)^{\frac{3}{p}}R^2(T)+\mathbb{F}_p(T)\rk]. 
 	\EEQS
 	
 	Letting $p \longrightarrow \infty$ completes the proof of Theorem \ref{thm2.8}.
 \end{proof}

 \section{Estimates on the spectral Energy function and Inertial Ranges}
 \label{sec: energyspectra}

 In this section we   study  the bounds of our spectral energy function $E(k,t)$ defined by \eqref{spectrafun}. 
 The results are presented in three theorems. 
 The first theorem ensures that the spectral energy remains bounded when the initial conditions and the non-homogeneous external forces satisfy certain conditions. 
 The second theorem estimates the time average   of the spectral energy;   it is shown that the average is always bounded and decays over time. 
 And finally, the third theorem gives the inertial range  bounds and formulates the conditions expected from the parameters, such as the dissipation rate, the universal constant and   viscosity coefficients so that the spectral energy decays accordingly with K-41. 
 This is done by comparing   $E(k,t)$ with Kolmogorov's spectral function $ E_K(k) $ given  by,
 
 \DEQSZ
 E_K(k) = C_0\epsilon^{\frac23}k^{ {-\frac{5}3}}, \label{k41likemhd}
 \EEQSZ
 
 \noindent defined over a range of wave numbers called the inertial range; where $C_0$ is a universal constant called Kolmogorov constant and $\epsilon$ is  {the} energy  dissipation rate. 
 
  {\begin{remark}\label{rem:k41likemhd}
 		Equation \eqref{k41likemhd} is similar to Equation (106) of \cite[pp.~ 267]{verma2004statistical} where $C_0$ and $\epsilon$ were referred to as Kolmogorov constant for MHD turbulence and energy flux respectively instead of Kolmogorov's constant and energy dissipation rate.
 \end{remark}}

 \del{ As discussed in section \ref{introduction}, K-41 theory   
 	postulates that the spectral energy of a turbulent flow  
 	decays proportional to
 	
 	\DEQSZ
 	C_0\epsilon^{\red{2/3}}k^{-5/3} \label{enrgysfl}.
 	\EEQSZ
 	
 	Even if  the K-41 theory was established for non magnetized fluid flows, several studies indicated that MHD flows also satisfy K-41 over a certain range, \cite{verma2004statistical,verma1999mean,biskamp1994cascade,beresnyak2011spectral,marsch1991turbulence}.
 	
 	Of great scientific interest is the question of rigorous mathematical proof of \red{Kolmogorov's} hypothesis under physically acceptable conditions. 
 	Biryuk and Craig in \cite{biryuk2012bounds} developed a rigorous mathematical technique to establish Kolmogorov's spectral theory for Navier-Stokes equations with certain smallness and regularity conditions imposed on the   data. 

 	Therefore, our main task in this  section  will be deriving the K-41 theory from a purely mathematical point of view; the method  of Biryuk and Craig in \cite{biryuk2012bounds} is used for this  purpose.}

 Recall that  the spectral energy function for  {the MHD system} \eqref{mhdmain2} 
 
 \DEQSZ  
 \lk\lbrace \begin{array}{l l}
 	\partial_t u +(u\cdot\nabla)u+\nabla\pi - (b\cdot \nabla)b-\nu\Delta u = f_1 & (0,\infty)\times D, \\
 	\partial_t b +(u\cdot\nabla)b - (b\cdot \nabla)u-\eta\Delta b = f_2 & (0,\infty)\times D ,\\
 	\text{div } u = {\rm div } b = 0 & D, \\
 	u|_{t=0} = u_0,  \quad b|_{t=0}=b_0  & D,
 \end{array} \rk.   \label{mhdmain2}
 \EEQSZ  
 
 \noindent is  given by the spherical integral 
 
 \DEQSZ
 E(k,t)=\int_{|\xi|=k}(|\widehat{u}(\xi,t)|^2+|\widehat{b}(\xi,t)|^2){\rm d} S(\xi) \label{alt44},
 \EEQSZ
 
 \noindent where $0\leq k<\infty$ is a radial coordinate in Fourier space.

 \begin{Theorem}\label{thm3.1}
 	Let the assumptions of Theorem \ref{thm2.7} hold,  $f_i\equiv 0$  for all $i=1,\,2$ and the initial data $(u_0,  b_0)\in  B_R(0),$ where $R$   satisfies  \eqref{eq: formainthm}. 
 	Then, the estimate
 	\DEQSZ
 	E(k,t)\leq 4\pi R_1^2, \label{49}
 	\EEQSZ
 	holds for all $k$ and all $t,$ where $R_1$ is as in Theorem \ref{thm2.7}.
 	Moreover,
 	when $f_i\not\equiv0$ for some $i=1,2$, \eqref{49} still holds with $R_1$ replaced by $R_1(t)$ which is still finite and possibly grows in time.  
 \end{Theorem}
 
 \begin{proof}[Proof of Theorem \ref{thm3.1}]
 	When $f_i\equiv 0,$ we have from   \eqref{alt44} and Theorem \ref{thm2.7} that 
 	\DEQS
 	E(k,t) &=& \int\limits_{|\xi|=k}(|\widehat{u}(\xi,t)|^2+|\widehat{b}(\xi,t)|^2)dS(\xi)\\
 	&\leq & \int\limits_{|\xi|=k}\frac{R_1^2}{k^2}dS(\xi) \\
 	&=& 4\pi R_1^2.
 	\EEQS
 	
 	\noindent	Here we used the fact that the surface area of a sphere with radius $k$ is equal to $4\pi k^2.$
 	
 	 {When}  the external forces on the system, $f_i\not\equiv0, $ for some $i=1,2$ the proof above remains same with  $R_1$   replaced with $R_1(t)$. 
 	With this we  complete the proof.
 \end{proof}
 \begin{remark}
 	Theorem \ref{thm3.1} implies that when no external force is applied to the system, the spectral energy remains uniformly bounded through out the entire process.
 \end{remark}
 \begin{Theorem}\label{thm3.2}
 	Suppose the initial data $(u_0,\, b_0)\in B_R(0),$ where $ R$ satisfy the conditions of Theorem \ref{thm3.1}
 	and the forces $f_i\in L_{loc}^\infty ([0,\infty]; H^{-1}(D)\cap L^2(D))$ for $i=1,2$
 	is bounded as it appears  in  \eqref{29}.
 	Then for every $T$, we have 
 	
 	\DEQSZ
 	\frac{1}{T}\int\limits_0^TE(k,t)dt\leq \frac{4\pi R_2^2(T)}{\min(\nu,\eta) Tk^2}\label{51},
 	\EEQSZ
 	where $R_2$ is as in Theorem \ref{thm2.8}.
 \end{Theorem}
 
 \begin{proof}[Proof of Theorem \ref{thm3.2}]
 	
 	From \eqref{alt44} we have,
 	
 	\DEQS
 	\frac{1}{T}\int\limits_0^TE(k,t)dt &=& \frac{1}{T}\int\limits_0^T\int\limits_{|\xi|=k}(|\widehat{u}(\xi,t)|^2+|\widehat{b}(\xi,t)|^2)dS(\xi)dt\\
 	&=& \frac{1}{T}\int\limits_0^T\int\limits_{|\xi|=k}(|\widehat{\chi}_k(\xi)\widehat{u}(\xi,t)|^2+|\widehat{\chi}_k(\xi)\widehat{b}(\xi,t)|^2)dS(\xi)dt\\
 	& \leq& \int\limits_{|\xi|=k}\frac{1}{T}\frac{ R_2^2(T)}{\min(\nu,\eta) k^4}dS(\xi) \quad \\
 	& \leq & 4\pi k^2 \frac{ R_2^2(T)}{\min(\nu,\eta) T k^4}  = \frac{ 4\pi R_2^2(T)}{\min(\nu,\eta) T k^2}.
 	\EEQS
 	
 	\noindent Here we used \eqref{28} and the proof is complete.
 \end{proof}
 
 %
 %
 
 \begin{Theorem}\label{smashthm} 
 	Let the assumptions of Theorem \ref{thm3.1} and Theorem \ref{thm3.2} hold. 
 	Then the following  
 	are true about  the inertial range of \eqref{mhdmain}:
 	\begin{enumerate}
 		\item  {Inequality \eqref{53} is a necessary condition on the parameters so that $E(k,t)$ exhibits  K-41 like phenomenon.}
 		\DEQSZ
 		\left(\min(\nu,\eta)\right)^{5/6}C_0\epsilon^{2/3 }\leq 4\pi  \left(\frac{R_2(T)}{\sqrt{T}} \right)^{5/3}R_1^{\frac{1}{3}}(T)\label{53}.
 		\EEQSZ
 		
 		\item An absolute lower bound for the inertial range is given by
 		
 		\DEQSZ
 		{k}_1= \frac{C_0^{3/5}\epsilon^{2/5}}{(4\pi R_1^2)^{3/5}}\label{55}.
 		\EEQSZ
 		
 		\item An absolute upper bound for the inertial range is given by
 		
 		\DEQSZ
 		{k}_2=\left(\frac{4\pi}{C_0\min(\nu,\eta)}\right)^3\frac{1}{\epsilon^2}\frac{R_2^6(T)}{T^3} \label{57}.
 		\EEQSZ
 		
 	\end{enumerate}
 \end{Theorem}
 \begin{proof}[Proof of Theorem \ref{smashthm}]
 	 {Define set $S$ by
 		
 		\DEQSZ
 		S:=\left\lbrace (k,E):0\leq E(k,\cdot)\leq 4\pi R_1^2,\,  0\leq E(k,\cdot)\leq \frac{ 4\pi R_2^2}{\min(\nu,\eta) Tk^2} \right\rbrace.
 		\EEQSZ
 		
 	}
 	
 	Let 
 	
 	\DEQS
 	A:= \{(k,E): E=E_K(k)\}\cap S,
 	\EEQS   
 	
 	\noindent be part of the graph of $E_K(k)$ that lies in region $S$. 
 	Figure \ref{fig1} shows how sets $S$ and $A$ are related.
 	\begin{figure}[h!]
 		\centering
 		\includegraphics[scale=0.4]{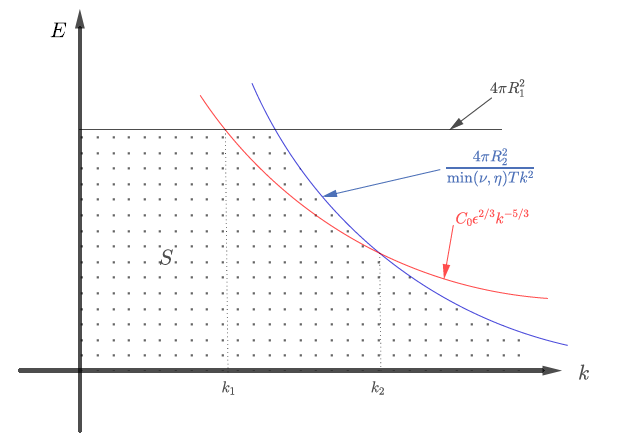}
 		\caption{Sketch of region $S$ and set $A$.}\label{fig1}
 	\end{figure}
 	
 	Due to Theorem \ref{thm3.1} we know that the spectral energy of our system is bounded from above by   $ 4\pi R_1^2$ when $f_i\equiv 0$ for all $i=1,2$ or $ 4\pi R_1^2(T)$ $f_i\not\equiv 0$ for some $i=1,2$. Furthermore, from  Theorem \ref{thm3.2}  the time average is bounded  by 
 	
 	\DEQS
 	\frac{4\pi R_2^2(T)}{\min{(\nu,\eta)}Tk^2}.
 	\EEQS
 	
 	Thus, set $S$ represents the behavior of  the   function $E(k,t)$ and set $A$ is a set where $E(k,t)$ behaves accordingly with K-41. 
 	Therefore if $A=\emptyset$ then $E(k,t)$ does not exhibit K-41 like phenomenon.
 	
 	Note that for $A$ to be non-empty  the point where graphs of  $E_K(k)$ and $\frac{ 4\pi R_2^2}{\min(\nu,\eta) Tk^2}$ must intersect   below the line  $E=4\pi R_1^2$, as  in Figure \ref{fig1}, and the
 	intersection occurs when
 	
 	\DEQS
 	C_0\epsilon^{2/3 } k^{-5/3} &=& \frac{ 4\pi R_2^2(T)}{\min(\nu,\eta) Tk^2} \\
 	\implies k &=& \left( \frac{4\pi R_2^2(T)}{\min(\nu,\eta) T C_0 \epsilon^{2/3}} \right)^3.
 	\EEQS
 	
 	%
 	
 	Moreover, the graph of $ E_K(k)$ intersects the line $E=4\pi R_1^2$ below the graph of $\frac{ 4\pi R_2^2}{\min(\nu,\eta) Tk^2}$, as   in Figure \ref{fig1}, which occurs when
 	
 	\DEQS
 	C_0\epsilon^{2/3 } k^{-5/3} &=& 4\pi R_1^2 \\
 	\implies  k &=& \left(\frac{ 4\pi R_1^2}{C_0\epsilon^{\frac{2}{3} } } \right)^{\frac{-3}{5}} .
 	\EEQS
 	
 	Therefore, $E_K(k)$ enters region $S$ at $k = \left(\frac{ 4\pi R_1^2}{C_0\epsilon^{\frac{2}{3} } } \right)^{\frac{-3}{5}} $ and leaves at $k=  \left( \frac{4\pi R_2^2(T)}{\min(\nu,\eta) T C_0 \epsilon^{\frac{2}{3}}} \right)^3$.
 	
 	Now set, 
 	
 	\DEQS
 	{ k}_1 = \left(\frac{ 4\pi R_1^2}{C_0\epsilon^{2/3 } } \right)^{-3/5},\qquad
 	{ k}_2=  \left( \frac{4\pi R_2^2(T)}{\min(\nu,\eta) T C_0 \epsilon^{2/3}} \right)^3 .
 	\EEQS
 	
 	Thus the portion of the graph of $E_K(k)$ remains in region $S$ as long as $k$ is between $ {k}_1$ and $ {k}_2$ and $k_1\leq k_2$. 
 	Hence, $A$ remains non-empty only when $k\in[k_1,k_2]$. 
 	
 	Observe from Figure \ref{fig1} that if we push the graph of $ \frac{ 4\pi R_2^2(T)}{\min(\nu,\eta) Tk^2} $ to the left so that it intersects $E_K(k)$ above the graph of $4\pi R_1^2$, then we get $k_1 >  k_2$ and the graph of $E_K(k)$ will not pass through region $S$ which in turn gives $A=\emptyset$.
 	
 	Therefore, for the flow model \eqref{mhdmain2} exhibit K-41 like MHD phenomenon we need the necessary condition 
 	
 	\DEQS
 	\left(\frac{ 4\pi R_1^2}{C_0\epsilon^{2/3 } } \right)^{-3/5} &\leq  \left( \frac{4\pi R_2^2(T)}{\min(\nu,\eta) T C_0 \epsilon^{2/3}} \right)^3 ,
 	\EEQS
 	
 	to be satisfied. 
 	Hence, 
 	
 	\DEQS
 	C_0 \min(\nu,\eta)^{5/6}\epsilon^{2/3 }&
 	\leq& 4\pi \left(\frac{R_2(T)}{\sqrt{T}} \right)^{5/3}R_1^{\frac{1}{3}}(T).
 	\EEQS
 	
 	This completes the proof Theorem \ref{smashthm}.  
 \end{proof}

 \section{Conclusion}
 
 In this work we have investigated the Leray weak solution of the deterministic MHD model \eqref{mhdmain} for the K-41 like MHD phenomenon in the presence and absence of nonhomogeneous external forces.
 In the process we have shown in Section \ref{fR3} that when the external forces $f_1$ and $f_2$ are identically $ 0 $ the solution field $(u,b)$ is bounded in the Fourier space (Theorem \ref{thm2.7} and Theorem \ref{thm2.8}).
 We have also shown that the spectral energy of the system $E(k,t)$ is bounded. In fact,  when $f_i\equiv 0$ for $i=1,2$ the bound is uniform  (Theorem \ref{thm3.1}) and the   average in time is also bounded, the bound decreases in time and decays    proportional to $k^{-2}.$   
 When $f_i\not\equiv 0$ for some $i=1,2$ the bonds of $E(k,t)$ possibly depend on time. 
 The other important result of this work is explicit formulation of the inertial range bounds and setting the necessary condition on the parameters for the model to behave accordingly with K-41 (Theorem \ref{smashthm}). The lower bound 
 \DEQS 
 { k}_1 = \left(\frac{ 4\pi R_1^2}{C_0\epsilon^{2/3 } } \right)^{-3/5},
 \EEQS 
 is a constant in time when $f_i\equiv0$ for $i=1,2$ and possibly decreases in time when $f_i\not\equiv0$ for some $i=1,2$. 
 The upper bound of the inertial range 
 $${ k}_2=  \left( \frac{4\pi R_2^2(T)}{\min(\nu,\eta) T C_0 \epsilon^{2/3}} \right)^3, $$
 decreases in time when $f_i\equiv 0$ for $i=1,2$ and will remain decreasing as long as the $R_2\propto T^\alpha$ and $\alpha<1/2.$
 For the case where  $f_i \equiv 0$ for $i=1,2$, $R_1$ and $R_2$ are constants independent of time and at time $T=T_0$, where 
 
 \DEQSZ
 T_0:=\frac{(4\pi)^{6/5}R_2^2 R_1^{2/5}}{\epsilon^{4/5}C_0^{6/5}\min(\nu,\eta) } \label{maxtime}
 \EEQSZ
 
 \noindent we get $ {k}_1 =  {k}_2$.
 This means that for any time $T> T_0$, the spectral range is empty. 
 
 Consequently, time $T_0$ appears to be the maximal time to exhibit K-41 in the system.
 
 If we assume that the dissipation rate is time dependent then \eqref{maxtime} gives
 
 \DEQSZ
 \epsilon(T_0) = \frac{(4\pi)^{3/2}R_1^{1/2}R_2^{5/2}}{T_0^{5/4}\min(\nu,\eta)^{5/4}C_0^{3/2}}\label{epsilonmax}.
 \EEQSZ
 
 The time $T_0$ being the maximal time,  \eqref{epsilonmax} must be the minimum dissipation rate to maintain a spectral behavior.

 \section*{Acknowledgments}
 The author wishes to acknowledge Prof M Sango for suggesting the problem and for his highly insightful discussions. The author also wishes to thank the University of Pretoria where the paper was written and  of DST-NRF Center of Exelence in Mathematical and Statistical Sciences (CoE-MaSS) for their financial support during my stay at the University of Pretoria.


\begin{thebibliography}{999}
 	\bibitem{arsenio2014derivation}
 	Ars\'enio, D., Ibrahim, S., Masmoudi, N. A derivation of the magnetohydrodynamics stystem from Navier-Stokes-Maxwell systems. {\em Arch. Rational Mech. Anal} {\bf 2015}, {\em 216}, 767--812.
 	
 	\bibitem{alfonsi2009reynolds}
 	Alfonsi, G. Reynolds-averaged Navier--Stokes equations for turbulence modeling. {\em Applied Mechanics Reviews} {\bf 2009}, {\em 62}
 	
 	\bibitem{argyropoulos2015recent}
 	Argyropoulos, C.D; Markatos, N.C. Recent advances on the numerical modelling of turbulent flows. {\em Applied Mathematical Modelling} {\bf 2015}, {\em 39}, 693--732.
 	
 	\bibitem{beresnyak2012basic}
 	Beresnyak, A. Basic properties of magnetohydrodynamic turbulence in the inertial range. {\em Monthly Notices of the Royal Astronomical Society.} {\bf 2012}, {\em 422}, 3495--3502.
 	
 		\bibitem{boldyrev2005spectrum}
 	Boldyrev, S. On the spectrum of magnetohydrodynamic turbulence. {\em The Astrophysical Journal Letters} {\bf 2005}, {\em 626}, L37-L40
 	
 		\bibitem{beresnyak2011spectral} 	
 	Beresnyak, A.  Spectral slope and Kolmogorov constant of MHD turbulence. {\em Phys. Rev. Lett.} {\bf 2011}, {\em 106}, 075001.
 	
 		\bibitem{bahouri2011fourier}
 	Bahouri, H.; Chemin, J.-Y.; Danchin, R. {\em Fourier analysis and nonlinear partial differential equations.} Springer Science \& Business Media, 2011.
 	
 	\bibitem{biskamp1994cascade}
 	Biskamp, D. Cascade models for magnetohydrodynamic turbulence. {\em Phys. Rev. E.} {\bf 1994}, {\em 50}, 2702.
 	
 	\bibitem{biryuk2012bounds}
 	Biryuk, A.; Craig, W. Bounds on Kolmogorov spectra for the Navier--Stokes equations. {\em Phys. D.} {\bf 2012}, {\em 241}, 426--438.

 	\bibitem{chen2012kolmogorov}
 	Chen, G-Q; Glimm, J.  Kolmogorov's Theory of Turbulence and Inviscid Limit of the {N}avier-{S}tokes Equations in $\mathbb{R}^3$. {\em Commun. Math. Phys.} {\bf  2012}, {\em 310}, 267--283.
 	
 		\bibitem{Chandrasekhar1955}
 	Chandrasekhar, S. Hydromagnetic Turbulence. I. A Deductive Theory. {\em Proc. R. Soc. London, Ser. A.} {\bf 1955}, {\em 233}, 322--330.
 	
 		\bibitem{cannone2006cauchy}
 	Cannone, M; Miao, C. X.; Prioux, N.; Yuan, B. Q. The Cauchy problem for the magneto-hydrodynamic system. {\em Banach Center Publ.} {\bf 2006}, {\em 74}, 59.
 	
 	 	\bibitem{duchon2000inertial}
 	Duchon, J., Robert, R. Inertial energy dissipation for weak solutions of incompressible Euler and Navier-Stokes equations. {\em Nonlinearity} {\bf 2000}, {\em 13} 249--255.
 	
 	\bibitem{dobrowolny1980fully}
 	Dobrowolny, M.; Mangeney, A.; Veltri, P. Fully developed anisotropic hydromagnetic turbulence in interplanetary space. {\em Phys. Rev. Lett.} {\bf 1980}, {\em 45}, 144.
 	
 	\bibitem{davidson2011voyage}
 	Davidson, P. A.; Kaneda, Y.; Moffatt, K.; Sreenivasan, K. R. {\em A voyage through turbulence}; Cambridge University Press, 2011.
 	
 	\bibitem{eyink2001dissipation}
 	Eyink, G. L. Dissipation in turbulent solutions of 2D Euler equations. {\em Nonlinearity} {\bf 2001}, {\em 14}, 787--802.
 	
 		\bibitem{filho2006weak}
 	Filho, M. C. L., Mazzucato, A.L., Lopes, H. J. N. Weak solutions, renormalized solutions and enstrophy defects in 2D turbulence. {\em Arch. Rational Mech. Anal} {\bf 2006}, {\em 179}, 353--387.
 	
 	\bibitem{goedbloed1998derivation}
 	Goedbloed, J.P. `Derivation'of the MHD equation. {\em Fusion Technology} {\bf 1998}, {\em 33}, 97--104.
 	
 	\bibitem{hormander1983analysis}
 	H{\"o}rmander, L. {\em The Analysis of Linear Partial Differential Operators: Vol.: 1.: Distribution Theory and Fourier Analysis.} Springer-Verlag, 1983.
 	
 		\bibitem{iroshnikov1964turbulence}
 	Iroshnikov, P. S. Turbulence of a conducting fluid in a strong magnetic field. {\em Soviet Astronomy.} {\bf 1964}, {\em 7}, 566.
 	
 		\bibitem{jackson2007osborne}
 	Jackson, D.; Launder, B. Osborne {R}eynolds and the publication of his papers on turbulent flow. {\em Annu. Rev. Fluid Mech.} {\bf 2007}, {\em 39}, 19--35.
 	
 		\bibitem{kolmogorov1941degeneration}
 	Kolmogorov, A. N. On the degeneration of isotropic turbulence in an incompressible viscous fluid. {\em Dokl. Akad. Nauk SSSR.} {\bf 1941}, {\em 31}, 319--323.
 	\bibitem{kolmogorov1941dissipation}
 	Kolmogorov, Andrey Nikolaevich. Dissipation of energy in locally isotropic turbulence. {\em Dokl. Akad. Nauk SSSR.} {\bf 1941}, {\em 32}, 16.
 	\bibitem{kolmogorov1941equations}
 	Kolmogorov, A. N. Equations of turbulent motion in an incompressible fluid. {\em Cr Acad. Sci. URSS.} {\bf 1941}, {\em 30}, 301--305.
 	\bibitem{kolmogorov1941local}
 	Kolmogorov, A. N. The local Structure of turbulence in incompressible viscous fluid for very large {R}eynolds numbers. {\em Dokl. Akad. Nauk SSSR.} {\bf 1941}, {\em 30}, 299--303
 	
 	\bibitem{kraichnan1965inertial}
 	Kraichnan, R. H. Inertial-range spectrum of hydromagnetic turbulence. {\em Phys. Fluids (1958-1988).} {\bf 1965}, {\em 8}, 1385--1387.
 	\bibitem{kraichnan1971inertial}
 	Kraichnan, R. H. Inertial-range transfer in two-and three-dimensional turbulence. {\em J. Fluid Mech.} {\bf 1971}, {\em 47}, 525--535.
 	
 	 	\bibitem{kato1984strongl}
 	Kato, T. StrongL p-solutions of the Navier-Stokes equation inR m, with applications to weak solutions. {\em Mathematische Zeitschrift} {\bf 1984}, {\em 187}, 471--480.
 	
 	 	\bibitem{leray1931systeme}
 	Leray, J. Sur le systeme d'{\'e}quations aux d{\'e}riv{\'e}s partielles qui r{\'e}git l'{\'e}coulement premanant des fluids visqueux. {\em C. R. Math. Acad. Sci. Paris.} {\bf 1931}, {\em 192}, 1180--1182.
 	\bibitem{leray1933etude}
 	Leray, J. Etude de diverses {\'e}quations int{\'e}grales non lin{\'e}aires et de quelques probl{\'e}mes que pose l'hydrodynamique. {\em Th{\'e}ses fran{\c{c}}aises de l'entre-deux-guerres.} {\bf 1933}, {\em 142}, 1--82.
 	
 		\bibitem{millionshchikov1939decay}
 	Millionshchikov, M. D. Decay of homogeneous isotropic turbulence in viscous incompressible fluids. {\em Dokl. Akad. Nauk SSSR.} {\bf 1939}, {\em 22}, 236.
 	
 	 	\bibitem{moninyaglomI}
 	Monin, A. S.; Yaglom, A. M. {\em Statistical Fluid Mechanics: Mechanics of Turbulence}; John L. L., Eds.;  Dover Publications, Inc., Mineola, New York, 1975.
 	\bibitem{monin2013statisticalII}
 	Monin, A. S.; Yaglom, A. M. {\em Statistical fluid mechanics, volume II: mechanics of turbulence}. Dover Publications, 2007.
 	
 		\bibitem{marsch1991turbulence}
 	Marsch, E. Turbulence in the solar wind. In {\em Reviews in Modern Astronomy.} Springer, 1991, 145--156.
 	
 	\bibitem{marsch1990radial}
 	Marsch, E.; Tu, C-Y. On the radial evolution of MHD turbulence in the inner heliosphere. {\em Journal of Geophysical Research: Space Physics} {\bf 1990}, {\em 95}, 8211--8229.
 	
 		\bibitem{matthaeus1989extended}
 	Matthaeus, W.H.,   Zhou, Y.   Extended inertial range phenomenology of magnetohydrodynamic turbulence. {\em Physics of Fluids B: Plasma Physics (1989-1993).} {\bf 1989}, {\em 1}, 1929 -- 1931
 	
 	 	\bibitem{peng2008rigorous}
 	Peng, Y.-J., Wang, S. Rigourous dervation of incompressible e-MHD equations from compressible Euler-Maxwell equations. {\em SIAM J. Math. Anal.} {\bf 2008}, {\em 40},  540--565. 
 	
 		\bibitem{Prandtl1904}
 	Prandtl, L. On the motion of fluids with very little friction. In {\em Erly developments of modern aerodynamics}; Amer Inst of Aeronautics; 1904, 77--87.
 	
 	
 
 	\bibitem{reynolds1883experimental}
 	Reynolds, O. An experimental investigation of the circumstances which determine whether the motion of water shall be direct or sinuous, and of the law of resistance in parallel channels. {\em Proc. R. Soc. London} {\bf 1883}, {\em 35}, 84--99.
 	
 	\bibitem{reynolds1894dynamical}
 	Reynolds, O. On the dynamical theory of incompressible viscous fluids and the determination of the criterion. {\em Proc. R. Soc. London}, {\bf 1894},  {\em 56}, 40--45.
 	
 	 	\bibitem{serrin1964local}
 	Serrin, J. Local behavior of solutions of quasi-linear equations. {\em Acta Math.} {\bf 1964}, {\em 111}, 247--302.
 	
 		\bibitem{taylor1935statistical}
 	Taylor, G. I. Statistical theory of turbulence.  {\em Proc. R. Soc. London, Ser. A.} {\bf 1935}, {\em 151}, 421--444.
 	\bibitem{taylor1921experiments}
 	Taylor, G. I. Experiments with rotating fluids. {\em  Proc. R. Soc. London, Ser. A.} {\bf 1921}, {\em 100}, 114--121.
 	
 	

 	\bibitem{tikhomirov1991selected}
 	Tikhomirov, V.M. {\em Selected Works of {A}{N} {K}olmogorov: Volume I: Mathematics and Mechanics.} Springer Science \& Business Media, 1991.
 	
 	
 		\bibitem{von1948progress}
 	Von Karman, T. Progress in the statistical theory of turbulence. {\em Proc. Natl. Acad. Sci. U. S. A.} {\bf 1948 }, {\em 34}, 530.
 	\bibitem{von1949concept}
 	Von Karman, T.; Lin, C.C. On the concept of similiarity in the theory of isotropic turbulence. {\em Rev. Mod. Phys.} {\bf 1949}, {\em 21}, 516.
 	
 	\bibitem{verma2004statistical}
 	Verma, M. K. Statistical theory of magnetohydrodynamic turbulence: recent results. {\em Phys. Rep.} {\bf 2004}, {\em 401}, 229--380.
 	\bibitem{verma1999mean}
 	Verma, M. K. Mean magnetic field renormalization and Kolmogorov's energy spectrum in magnetohydrodynamic turbulence. {\em Phys. Plasmas.} {\bf 1999}, {\em 6}, 1455--1460.
 	
 	 	\bibitem{wolff2003lectures}
 	Wolff, T. H; Laba, I.; Shubin, C. {\em Lectures on harmonic analysis.} American Mathematical Society, 2003.
 	
 	
 	
 	
 	
 	
 	
 	
 	
 	
 	
 	
 	
 	
 	
 	
%
 	
 	 

%
 	
 	
 	

 	
 	
 	
 	
%
 \end{thebibliography}
\end{document}